\documentclass[12pt,reqno,a4paper]{amsart}
\usepackage{amsmath,amssymb,amsthm}
\usepackage[margin=1in,nohead,includefoot,footskip=12pt]{geometry}
\usepackage{cite}
\usepackage{tikz}
\usepackage{hyperref}
\hypersetup{
colorlinks=true,
linkcolor=blue,
urlcolor=red,
citecolor=green
}

\allowdisplaybreaks[4]

\newtheorem{theorem}{Theorem}[section]
\newtheorem{lemma}[theorem]{Lemma}
\newtheorem{proposition}[theorem]{Proposition}
\newtheorem{corollary}[theorem]{Corollary}
\numberwithin{equation}{section}
\newtheorem{remark}[theorem]{Remark}

\begin{document}

\title{Uniqueness of positive solutions of the Lane--Emden equation in smooth convex planar domains}

\author{Hong-Ge Chen}
\address{[Hong-Ge Chen] School of Mathematics and Statistics, Key Laboratory of Nonlinear Analysis Applications (Ministry of Education), Central China Normal University, Wuhan 430079, P.R. China}
\email{hongge\_chen@whu.edu.cn}

\author{Peng Luo}
\address{[Peng Luo] School of Mathematics and Statistics, Key Laboratory of Nonlinear Analysis Applications (Ministry of Education), Central China Normal University, Wuhan 430079, P.R. China}
\email{pluo@ccnu.edu.cn}

\author{Shuangjie Peng}
\address{[Shuangjie Peng] School of Mathematics and Statistics, Key Laboratory of Nonlinear Analysis Applications (Ministry of Education), Central China Normal University, Wuhan 430079, P.R. China}
\email{sjpeng@mail.ccnu.edu.cn}

\author{Fulin Zhong}
\address{[Fulin Zhong] School of Mathematics and Statistics, Central China Normal University, Wuhan 430079, P.R. China}
\email{flzhong@mails.ccnu.edu.cn}

\author{Yang Zhou}
\address{[Yang Zhou] School of Mathematics and Statistics, Key Laboratory of Nonlinear Analysis Applications (Ministry of Education), Central China Normal University, Wuhan 430079, P.R. China}
\email{yzhou@ims.cuhk.edu.hk}

\begin{abstract}
We study positive solutions of the Lane--Emden problem
\begin{equation*}
\begin{cases}
-\Delta u=u^p &\text{in }\Omega,\\
u=0 & \text{on }\partial\Omega,
\end{cases}
\end{equation*}
where $p>1$ and $\Omega$ is a smooth bounded convex domain in $\mathbb{R}^2$. We prove that every positive solution is a ground state. In particular, the problem has a unique positive solution $u$, which is nondegenerate and has Morse index one. Moreover, the Hessian of $u^{-\frac{p-1}{2}}$ is positive definite in $\Omega$. This settles the planar case of Dancer's conjecture proposed in 1988.

\vspace{0.25cm}

\vspace{0.25cm}
\end{abstract}

\maketitle

\section{Introduction}
In this paper, we consider the uniqueness of positive solutions of the well known Lane--Emden problem
\begin{equation}\label{eq:le}
\begin{cases}
-\Delta u=u^p &\text{in }\Omega,\\
u=0 & \text{on }\partial\Omega.
\end{cases}
\end{equation}
where $p>1$ and $\Omega$ is a smooth bounded domain in $\mathbb{R}^2$.

\vskip 0.1cm

When $\Omega$ is a disk, the celebrated theorem of Gidas, Ni, and Nirenberg~\cite{gnn} implies that every positive solution is radially symmetric. Uniqueness then follows from a scaling argument for the radial equation. Later, in 1980s, Dancer \cite{dancer} proved uniqueness on planar domains symmetric and convex with respect to two orthogonal directions. He \cite{dancer,dancer89} also constructed contractible nonconvex domains with more than one positive solution. In the plane, Dancer's conjecture \cite{dancer} takes the following form.

\vskip 0.2cm
\noindent\textbf{Conjecture A.} Lane--Emden problem
\begin{equation*}
\begin{cases}
-\Delta u=u^p &\text{in }\Omega,\\
u=0 & \text{on }\partial\Omega,
\end{cases}
\end{equation*}
has a unique positive solution for $p>1$ whenever $\Omega$ is convex (see also \cite{dancer,lions,lwz}).

\vskip 0.1cm

For general exponent $p$, another wonderful work is by Damascelli, Grossi and Pacella \cite{dgp}, who gave another proof on planar domains symmetric and convex with respect to two orthogonal directions based on PDE's method. They \cite{dgp} also established a more fundamental result for general nonlinearities, explaining why a solution is uniquely determined by its value at the center.

\vskip 0.1cm

For general domains, Lin \cite{lin} proved uniqueness on smooth convex planar domains for $p>1$ sufficiently close to $1$. As $p\to1^+$, the positive solution normalized by its maximum converges to the first Dirichlet eigenfunction with the same normalization.

\vskip 0.1cm
Scenario $p\to\infty$ has also been widely studied.
In the middle of 1990s, Ren and Wei \cite{ren94,ren96} proved that, along a subsequence as $p\to\infty$, least-energy solutions concentrate at a single interior point. Adimurthi and Grossi \cite{ag} identified the limiting peak height as $\sqrt{\mathrm{e}}$. For arbitrary positive solutions $u_p$, De Marchis, Ianni and Pacella \cite{dip} described concentration at finitely many interior points under the energy bound
\begin{equation*}
p\int_\Omega |\nabla u_p|^2\mathrm{d}x\leq C.
\end{equation*}
Recently, Kamburov and Sirakov \cite{ks} established this bound on smooth strictly star-shaped planar domains, in particular on smooth convex domains. On convex domains, a result of Grossi and Takahashi \cite{gtaka} excludes critical configurations of two or more distinct points for the associated interaction function. De Marchis, Grossi, Ianni and Pacella \cite{dm} used this result to exclude multiple concentration points. They computed the Morse index of single-peak solutions and proved that, for sufficiently large $p$, the positive solution on a smooth bounded convex planar domain is unique, nondegenerate and has Morse index one. Grossi, Ianni, Luo and Yan \cite{gily} also proved local uniqueness for single-peak solutions concentrating at the same nondegenerate critical point of the Robin function. This gives another proof of uniqueness on convex domains for large $p$.

\vskip 0.1cm
For every $p>1$, Lin \cite{lin} showed that $v^{-\frac{p-1}{2}}$ is convex for every positive ground state $v$ on a bounded convex planar domain and proved uniqueness of the positive ground state when the domain is smooth and strictly convex. Brasco and Franzina \cite{bf} proved uniqueness on every bounded convex planar domain without either assumption. It remained to determine whether a smooth bounded convex planar domain could also admit a positive solution of higher energy. In this paper, we prove that every positive solution is a ground state.

\vskip 0.1cm
To state our results, we introduce some notations. Denote the best Sobolev constant by
\begin{equation}\label{eq:sob}
S_p(\Omega)=\inf_{0\neq w\in H^1_0(\Omega)}
\frac{\displaystyle\int_\Omega |\nabla w|^2\mathrm{d}x}{\left(\displaystyle\int_\Omega |w|^{p+1}\mathrm{d}x\right)^{\frac{2}{p+1}}}.
\end{equation}
By compactness of the embedding $H^1_0(\Omega)\hookrightarrow L^{p+1}(\Omega)$, the infimum is attained. A nonnegative minimizer can be rescaled to solve problem \eqref{eq:le}. We call a positive solution attaining $S_p(\Omega)$ a ground state. For a positive solution $v$, let
\begin{equation*}
L_v=-\Delta-pv^{p-1}\quad\text{on }H^2(\Omega)\cap H^1_0(\Omega),
\end{equation*}
and denote by $\lambda_j(L_v)$ its eigenvalues, arranged in increasing order and counted with multiplicity. Then our first result is as follows.

\vskip 0.1cm
\begin{theorem}\label{thm:main}
Let $\Omega\subset\mathbb{R}^2$ be a bounded smooth convex domain and let $p>1$. Then every positive solution $u$ of problem \eqref{eq:le} satisfies
\begin{equation*}
\int_\Omega u^{p+1}\mathrm{d}x=\big(S_p(\Omega)\big)^{\frac{p+1}{p-1}}.
\end{equation*}
In particular, problem \eqref{eq:le} admits a unique positive solution $u$. Moreover,
\begin{equation*}
\lambda_1(L_u)<0<\lambda_2(L_u).
\end{equation*}
\end{theorem}

\begin{remark}
Theorem \ref{thm:main} settles the planar case of Dancer's conjecture (Conjecture A) in the affirmative. The conjecture remains open in dimensions $N \geq 3$.
\end{remark}

Now, using Theorem \ref{thm:main}, we can get following result on the qualitative properties on the solutions of problem \eqref{eq:le}.

\begin{corollary}\label{cor:conv}
Let $\Omega\subset\mathbb{R}^2$ be a bounded smooth convex domain and let $p>1$. Then every positive solution $u$ of problem \eqref{eq:le} satisfies
\begin{equation*}
D^2\left(u^{-\frac{p-1}{2}}\right)>0\quad\text{in }\Omega.
\end{equation*}
\end{corollary}

We outline the proof of Theorem \ref{thm:main}. Let $u$ be a positive solution of problem \eqref{eq:le} in $\Omega$, and assume first that $\partial\Omega$ has positive curvature $\kappa$. For $n(\theta)=(\cos\theta,\sin\theta)$, let $x(\theta)$ be the boundary point with outward unit normal $n(\theta)$. Since $\mathrm{d}\theta=\kappa\mathrm{d}s$, the Pohozaev identity for translations becomes
\begin{equation*}
\int_0^{2\pi}k_u(\theta)n(\theta)\mathrm{d}\theta=0,\quad
k_u(\theta):=\frac{|\partial_\nu u(x(\theta))|^2}{2\kappa(x(\theta))}.
\end{equation*}
Let $v$ be the positive ground state of problem \eqref{eq:le} in another bounded convex domain $\Omega'$ with positive boundary curvature. Using
$D^2(v^{-\frac{p-1}{2}})>0$ in $\Omega'$ from \cite{lin,gms}, we can translate and rescale $v$ to match the value and gradient of $u$ at each point. As in Section 3 of \cite{mira}, we compare $D^2u(x)$ with the Hessian of the translated and rescaled solution at $x$. If $k_u=k_v$, the boundary identities and Lemma \ref{lem:tensor} show that the two Hessians are equal. At points with the same outward normal, the normal derivatives are equal by the boundary identity, and the curvatures are equal since $k_u=k_v$. Thus the domains differ by a translation $c\in\mathbb{R}^2$, and the unique continuation principle yields
\begin{equation*}
\Omega=c+\Omega',\quad u(x)=v(x-c)\quad\text{for all }x\in\Omega.
\end{equation*}

Next, we choose a bounded convex domain $\Omega'$ with positive boundary curvature so that the positive ground state $v$ of \eqref{eq:le} in $\Omega'$ satisfies $k_v=k_u$. The kernel of the linearized map from the support function of $\Omega'$ to $k_v$ consists of the translation directions $a\cdot n$, $a\in\mathbb{R}^2$. After normalizing the support function, Lemma \ref{lem:inv} allows us to solve $k_v=k$ locally for functions $k$ satisfying
\begin{equation*}
\int_0^{2\pi}k(\theta)n(\theta)\mathrm{d}\theta=0.
\end{equation*}
Set $k_0=\max_{\mathbb{S}^1}k_u$. On a disk, the positive ground state $v$ of \eqref{eq:le} is radial, so $k_v$ is constant. We choose the radius so that $k_v\equiv k_0$. Starting from this disk, we prescribe
\begin{equation*}
k_t=(1-t)k_0+tk_u,\quad 0\leq t\leq1.
\end{equation*}
Let $[0,T)$ be the maximal interval on which the family can be continued from the disk. For $0\leq t<T$, let $h_t$ denote the normalized support function. Since $h_t>0$ and $k_u\leq k_0$, the first variation and the scaling law for the ground-state energy yield
\begin{equation*}
\frac{p+3}{p-1}\frac{\mathrm{d}}{\mathrm{d}t}\mathcal{E}(h_t)
=\int_0^{2\pi}h_t(\theta)(k_u(\theta)-k_0)\mathrm{d}\theta\leq0.
\end{equation*}
By Lemma \ref{lem:bound}, there exists a sequence $t_j\to T$ with $t_j<T$, such that $h_{t_j}$ converges to the support function $h_*$ of a domain with positive boundary curvature. By continuity of the boundary map, the positive ground state $v$ of \eqref{eq:le} in $\Omega_{h_*}$ satisfies $k_v=k_T$. If $T<1$, Lemma \ref{lem:inv} extends the family beyond $T$, a contradiction. Hence $T=1$ and $k_v=k_u$. Proposition \ref{prop:comp} then shows that $u$ is a translate of $v$.

\vskip 0.1cm
When the curvature vanishes somewhere on $\partial\Omega$, we perturb $\Omega$ to obtain a bounded convex domain with positive boundary curvature. If $\ker L_u\neq\{0\}$, we also vary the boundary on an arc where $\kappa>0$. Using the unique continuation principle, we can choose the normal variations on this arc so that the right-hand side of the linearized equation is orthogonal to $\ker L_u$. By the implicit function theorem, we obtain positive solutions $u_\varepsilon$ of \eqref{eq:le} in the perturbed domains $\Omega_\varepsilon=\Phi_\varepsilon(\Omega)$ such that
\begin{equation*}
u_\varepsilon\circ\Phi_\varepsilon\to u\quad\text{in }C^2(\overline{\Omega}).
\end{equation*}
Each $u_\varepsilon$ is a ground state by Proposition \ref{prop:pos}. Changing variables in the ground-state identity and using $S_p(\Omega_\varepsilon)\to S_p(\Omega)$, we obtain
\begin{equation*}
\int_\Omega u^{p+1}\mathrm{d}x=\big(S_p(\Omega)\big)^{\frac{p+1}{p-1}}.
\end{equation*}
Uniqueness and the spectral inequalities follow from Proposition \ref{prop:bf}.

\vskip 0.1cm

The paper is organized as follows. In Section \ref{sec:prelim}, we recall the properties of positive ground states and describe the boundary geometry of planar convex domains. In Section \ref{sec:comp}, we prove a rigidity result for domains with positive boundary curvature. In Section \ref{sec:real}, we construct a convex domain whose ground state has a prescribed boundary function. In Section \ref{sec:proof}, we complete the proof of the main results by allowing the boundary curvature to vanish.

\section{Ground states and boundary geometry}\label{sec:prelim}

In this section, we first recall the variational and spectral properties of positive ground states of problem \eqref{eq:le} on bounded convex domains in $\mathbb{R}^2$. For domains of class $C^{3,\delta}$, we derive the boundary identities for positive solutions. When the boundary curvature is positive, we also express the boundary point in terms of the support function.

\vskip 0.1cm

Fix $0<\delta<\min\{1,p-1\}$. Let $\Omega\subset\mathbb{R}^2$ be a bounded convex domain.

\begin{lemma}
The infimum in \eqref{eq:sob} is attained. Moreover, if
a nonnegative function $w$ is a minimizer of \eqref{eq:sob} with $\|w\|_{L^{p+1}(\Omega)}=1$, then $S_p(\Omega)^{\frac{1}{p-1}}w$ is a positive ground state of problem \eqref{eq:le}.
\end{lemma}
\begin{proof}
By the Sobolev inequality, $S_p(\Omega)>0$. Since the embedding $H^1_0(\Omega)\hookrightarrow L^{p+1}(\Omega)$ is compact, the infimum in \eqref{eq:sob} is attained by a nonnegative function $w$ with $\|w\|_{L^{p+1}(\Omega)}=1$. The Euler-Lagrange equation for $w$ is
\begin{equation*}
\begin{cases}
-\Delta w=S_p(\Omega)w^p &\text{in }\Omega,\\
w=0 &\text{on }\partial\Omega.
\end{cases}
\end{equation*}
In dimension two, $w\in H^1_0(\Omega)$ implies $w\in L^{2p}(\Omega)$, so $w^p\in L^2(\Omega)$. By the Dirichlet $H^2$ estimate on bounded convex domains (see Theorem 3.2.1.2 in \cite{Gr}) and the Sobolev embedding theorem, we have
\begin{equation*}
w\in H^2(\Omega),\quad w\in C^{0,\sigma}(\overline{\Omega})\quad\text{for every }0<\sigma<1.
\end{equation*}
Since $w^p$ is H\"older continuous, then using the Schauder estimate, we have $w\in C^2(\Omega)$. As $w\geq0$ and $w\not\equiv0$, the strong maximum principle gives $w>0$ in $\Omega$. Thus, $S_p(\Omega)^{\frac{1}{p-1}}w$ is a positive ground state of problem \eqref{eq:le}.
\end{proof}

\vskip 0.1cm

\begin{lemma}
Let $u\in C^2(\Omega)\cap C(\overline{\Omega})$ be a positive solution of problem \eqref{eq:le}.  Then $u$ is a ground state of problem \eqref{eq:le} if and only if
\begin{equation}\label{eq:least}
\int_\Omega u^{p+1}\mathrm{d}x=\big(S_p(\Omega)\big)^{\frac{p+1}{p-1}}.
\end{equation}
\end{lemma}
\begin{proof}

For $\varepsilon>0$, set $u_\varepsilon=(u-\varepsilon)_+$. The function $u_\varepsilon$ is locally Lipschitz in $\Omega$ and has compact support in $\Omega$. Hence $u_\varepsilon\in H^1_0(\Omega)$ and $\nabla u_\varepsilon=\mathbf{1}_{\{u>\varepsilon\}}\nabla u$ almost everywhere. Testing \eqref{eq:le} with $u_\varepsilon$, we obtain
\begin{equation*}
\int_{\{u>\varepsilon\}}|\nabla u|^2\mathrm{d}x=\int_\Omega u^p u_\varepsilon\mathrm{d}x.
\end{equation*}
USing the monotone convergence theorem, we have
\begin{equation}\label{eq:energy}
\int_\Omega|\nabla u|^2\mathrm{d}x=\int_\Omega u^{p+1}\mathrm{d}x<\infty.
\end{equation}
Since $0\leq u-u_\varepsilon\leq\varepsilon$ and $|\nabla u|^2\in L^1(\Omega)$, we get
\begin{equation*}
\|u-u_\varepsilon\|_{H^1(\Omega)}^2\leq\varepsilon^2|\Omega|+\int_{\{u\leq\varepsilon\}}|\nabla u|^2\mathrm{d}x\to0.
\end{equation*}
Hence $u\in H^1_0(\Omega)$. It follows from \eqref{eq:sob} and \eqref{eq:energy} that $u$ is a ground state of problem \eqref{eq:le} if and only if
\eqref{eq:least} holds.
\end{proof}

\begin{lemma}[See {\cite[Theorem 4.5]{bf}}]\label{luo-1}
Let $\Omega\subset\mathbb{R}^2$ be a bounded convex domain and let $p>1$. Then \eqref{eq:sob} has a unique positive minimizer $w$ satisfying $\|w\|_{L^{p+1}(\Omega)}=1$.
\end{lemma}

\begin{proposition}\label{prop:bf}
Let $\Omega\subset\mathbb{R}^2$ be a bounded convex domain and let $p>1$. Then problem \eqref{eq:le} has a unique positive ground state $v$, and
\begin{equation}\label{eq:spec}
\lambda_1(L_v)<0<\lambda_2(L_v).
\end{equation}
\end{proposition}

\begin{proof}
Let $v$ be a positive ground state of problem \eqref{eq:le} in $\Omega$ and set
\begin{equation*}
w=\big(S_p(\Omega)\big)^{-\frac{1}{p-1}}v.
\end{equation*}
Then $w$ minimizes \eqref{eq:sob}, and \eqref{eq:least} gives $\|w\|_{L^{p+1}(\Omega)}=1$. By Lemma \ref{luo-1}, $w$ is the unique positive minimizer with this normalization. Hence $v=\big(S_p(\Omega)\big)^{\frac{1}{p-1}}w$ is the unique positive ground state of problem  \eqref{eq:le} in $\Omega$.

\vskip 0.1cm

By the min-max principle, we have
\begin{equation*}
\lambda_1(L_v)=\min_{\substack{\phi\in H^1_0(\Omega)\\ \|\phi\|_{L^2(\Omega)}=1}}\int_\Omega \left(|\nabla\phi|^2-pv^{p-1}\phi^2\right)\mathrm{d}x.
\end{equation*}
Taking $\phi=\frac{v}{\|v\|_{L^2(\Omega)}}$ and using \eqref{eq:energy}, we obtain
\begin{equation*}
\lambda_1(L_v)\int_\Omega v^2\mathrm{d}x\leq\int_\Omega \left(|\nabla v|^2-pv^{p+1}\right)\mathrm{d}x=-(p-1)\int_\Omega v^{p+1}\mathrm{d}x<0.
\end{equation*}

Since $S_p(\Omega)w^{p-1}=v^{p-1}$, \cite[Lemma A]{bf} with $q=p+1$ gives $\lambda_2(L_v)\geq0$. Suppose that $\lambda_2(L_v)=0$, and let $\psi$ be a corresponding eigenfunction. For $X(x)=x-a$ with $a\in\mathbb{R}^2$, or for a constant vector field $X$, the boundary identity in the proof of Theorem 4.5 in \cite{bf} gives
\begin{equation*}
\int_{\partial\Omega}(X\cdot\nu)\partial_\nu w \partial_\nu\psi \mathrm{d}s=0.
\end{equation*}
By the convexity argument in that proof, we can choose $X$ so that $(X\cdot\nu)\partial_\nu\psi$ has a fixed sign almost everywhere on $\partial\Omega$ and is nonzero on a set of positive arclength. When the nodal line of $\psi$ meets $\partial\Omega$ at two points, the choice of $X$ uses the supporting lines at those points. Since $\partial_\nu w<0$ almost everywhere on $\partial\Omega$, the integral is nonzero. This contradiction proves $\lambda_2(L_v)>0$.
\end{proof}

Set $\beta=\frac{p-1}{2}$. We recall the following convexity result for positive ground states of problem \eqref{eq:le}.

\begin{proposition}[See {\cite[Theorem 2]{lin}}]\label{prop:lin}
Let $\Omega\subset\mathbb{R}^2$ be a bounded convex domain and let $p>1$. If $v$ is the positive ground state of problem \eqref{eq:le}, then
\begin{equation*}
D^2(v^{-\beta})\geq0\quad\text{in }\Omega.
\end{equation*}
\end{proposition}

\begin{lemma}\label{lem:conv}
Let $\Omega\subset\mathbb{R}^2$ be a bounded convex domain and let $p>1$. If $v$ is the positive ground state of problem \eqref{eq:le}, then
\begin{equation}\label{eq:conv}
D^2(v^{-\beta})>0\quad\text{in }\Omega.
\end{equation}
\end{lemma}
\begin{proof}
By Proposition \ref{prop:lin}, $G=v^{-\beta}$ is convex in $\Omega$. By $p=2\beta+1$ and $-\Delta v=v^p$, direct differentiation gives
\begin{equation*}
\begin{aligned}
\nabla G&=-\beta v^{-\beta-1}\nabla v,
\end{aligned}
\end{equation*}
and
\begin{equation*}
\begin{aligned}
\Delta G&=-\beta v^{-\beta-1}\Delta v+\beta(\beta+1)v^{-\beta-2}|\nabla v|^2\\
&=\beta v^\beta+\beta(\beta+1)v^{-\beta-2}|\nabla v|^2
=\frac{\beta+(1+\beta^{-1})|\nabla G|^2}{G}=:b(G,\nabla G).
\end{aligned}
\end{equation*}
Since $v>0$ in $\Omega$, the Schauder estimate \cite[Chapter 6]{gt} yields $v\in C^4(\Omega)$. For $\varphi(t)=t^{-\beta}$, we have $\varphi'<0<\varphi''$ on $(0,\infty)$. For each fixed $z$, $b(t,z)>0$ and $1/b(t,z)$ is affine in $t>0$. Since $G=\varphi(v)$ is convex, \cite[Corollary 2.7]{gms} yields $D^2G>0$ in $\Omega$.
\end{proof}

Let $v$ be a positive solution of problem \eqref{eq:le} in a bounded convex domain $\Omega$ of class $C^{3,\delta}$. Since $v^p$ is bounded, the $W^{2,r}$ estimate in \cite[Chapter 9]{gt} gives $v\in W^{2,r}(\Omega)$ for every $r>2$. Taking $r>\frac{2}{1-\delta}$ and applying the Sobolev embedding theorem, we obtain
\begin{equation*}
v\in C^{1,\delta}(\overline\Omega),\quad \nabla(v^p)=pv^{p-1}\nabla v\in C^{0,\delta}(\overline\Omega).
\end{equation*}
Let $\nu=(\nu_1,\nu_2)$ be the outward unit normal to $\partial\Omega$ and set $q_v=-\partial_\nu v$. We have $v^p\in C^{1,\delta}(\overline{\Omega})$, $\Delta v<0$ in $\Omega$ and $v=0$ on $\partial\Omega$. Using the Schauder estimate \cite[Chapter 6]{gt} and the Hopf lemma, we can get
\begin{equation}\label{eq:reg}
v\in C^{3,\delta}(\overline\Omega),\quad q_v>0\quad\text{on }\partial\Omega.
\end{equation}

For $i=1,2$, multiply $-\Delta v=v^p$ by $\partial_{x_i}v$ and integrate over $\Omega$. Since $v=0$ and $\nabla v=-q_v\nu$ on $\partial\Omega$, integration by parts gives
\begin{equation*}
\begin{aligned}
\int_\Omega v^p\partial_{x_i}v\mathrm{d}x&=\frac{1}{p+1}\int_{\partial\Omega}v^{p+1}\nu_i\mathrm{d}s=0,
\end{aligned}
\end{equation*}
and
\begin{equation*}
\begin{aligned}
\int_\Omega(-\Delta v)\partial_{x_i}v\mathrm{d}x&=\int_{\partial\Omega}\left(\frac{1}{2}|\nabla v|^2\nu_i-\partial_\nu v\partial_{x_i}v\right)\mathrm{d}s=-\frac{1}{2}\int_{\partial\Omega}q_v^2\nu_i\mathrm{d}s.
\end{aligned}
\end{equation*}
We obtain the Pohozaev identity for translations (see also \cite{pohoz,pucci})
\begin{equation}\label{eq:flux}
\int_{\partial\Omega}q_v^2\nu\mathrm{d}s=0.
\end{equation}

Let $\tau=(-\nu_2,\nu_1)$ be the unit tangent to $\partial\Omega$ and let $s$ denote arc length in the direction of $\tau$. With the curvature $\kappa:=\partial_\tau\nu\cdot\tau$, the Frenet formulas \cite[Section 1-5]{docarmo} take the form
\begin{equation} \label{eq:Fre}
\partial_\tau\nu=\kappa\tau,\quad \partial_\tau\tau=-\kappa\nu.
\end{equation}
Since $v=0$ on $\partial\Omega$, we have $\partial_\tau v=0$ and $\nabla v=(\partial_\nu v)\nu=-q_v\nu$ on $\partial\Omega$. Differentiating along $\tau$ and using \eqref{eq:Fre}, we obtain
\begin{equation*}
D^2v \tau=-(\partial_\tau q_v)\nu-q_v\partial_\tau\nu=-(\partial_\tau q_v)\nu-\kappa q_v\tau.
\end{equation*}
Since $v\in C^2(\overline{\Omega})$, the equation $\Delta v=-v^p$ extends to $\partial\Omega$, where $v=0$. Thus $\Delta v=0$ on $\partial\Omega$, and the boundary components of $D^2v$ satisfy
\begin{equation}\label{eq:hess}
\begin{aligned}
D^2v(\tau,\tau)&=(D^2v \tau)\cdot\tau=-\kappa q_v,\\
D^2v(\nu,\tau)&=(D^2v \tau)\cdot\nu=-\partial_\tau q_v,\\
D^2v(\nu,\nu)&=\Delta v-D^2v(\tau,\tau)=\kappa q_v.
\end{aligned}
\end{equation}

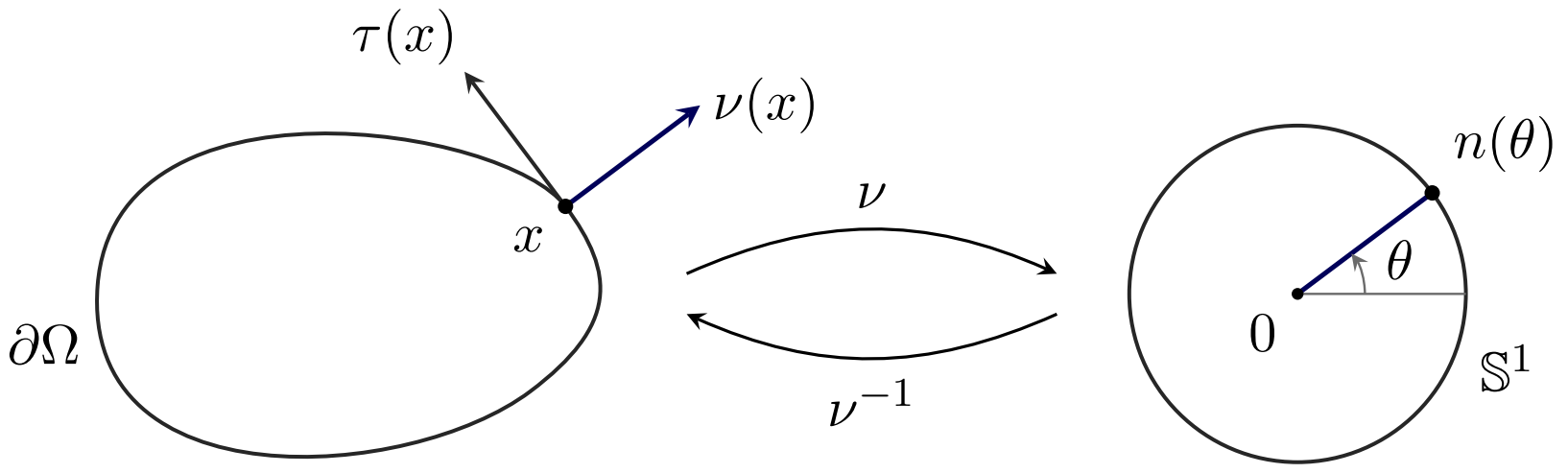
\begin{figure}[htbp]
\centering
\begin{tikzpicture}[
scale=1.16, >=stealth, font=\small,
outline/.style={draw=black!85, line width=0.75pt},
normal/.style={draw=blue!35!black, line width=0.95pt},
map/.style={->, line width=0.60pt},
guide/.style={draw=black!58, line width=0.45pt}
]
\coordinate (X) at (1.20,0.52);
\draw[outline] (X)
.. controls (0.84,1.00) and (-1.48,1.35) .. (-1.58,0.05)
.. controls (-1.68,-1.25) and (0.27,-1.12) .. (0.99,-0.58)
.. controls (1.51,-0.19) and (1.50,0.12) .. (X) -- cycle;
\node at (-1.90,-0.3) {$\partial \Omega$};
\draw[normal,->] (X) -- ++(0.8,0.6) node[right,inner sep=2pt] {$\nu(x)$};
\draw[outline,->] (X) -- ++(-0.6,0.8) node[above left,inner sep=1pt] {$\tau(x)$};
\fill (X) circle (1.3pt) node[below left] {$x$};

\draw[map] (1.92,0.12) to[bend left=24] node[midway,above] {$\nu$} (4.12,0.12);
\draw[map] (4.12,-0.12) to[bend left=24] node[midway,below] {$\nu^{-1}$} (1.92,-0.12);

\begin{scope}[shift={(5.55,0)}]
\pgfmathsetmacro{\ang}{atan(0.6/0.8)}
\draw[outline] (0,0) circle (1);
\node[right,inner sep=3pt] at (-24:1.08) {$\mathbb{S}^1$};
\draw[guide] (0,0) -- (1,0);
\draw[normal] (0,0) -- (0.8,0.6);
\draw[guide,->] (0.40,0) arc (0:\ang:0.40);
\node at ({\ang/2}:0.64) {$\theta$};
\fill (0,0) circle (1pt) node[below left] {$0$};
\fill (0.8,0.6) circle (1.3pt) node[above right] {$n(\theta)$};
\end{scope}
\end{tikzpicture}
\caption{The outward normal map.}
\label{fig:normal}
\end{figure}

Assume that $\kappa>0$ on $\partial\Omega$. Then the outward unit normal map $\nu:\partial\Omega\to\mathbb{S}^1$ is a $C^{2,\delta}$ diffeomorphism. As in Figure \ref{fig:normal}, set
\begin{equation*}
n(\theta)=(\cos\theta,\sin\theta),\quad x(\theta)=\nu^{-1}(n(\theta)).
\end{equation*}
Write $\tau(\theta)=\tau(x(\theta))=n'(\theta)$. Since $x'(\theta)$ is tangent to $\partial\Omega$, differentiating $\nu(x(\theta))=n(\theta)$ and using \eqref{eq:Fre}, we obtain $\kappa(x(\theta))x'(\theta)=\tau(\theta)$. Hence
\begin{equation} \label{eq:x'-theta}
x'(\theta)=\frac{\tau(\theta)}{\kappa(x(\theta))},\quad \frac{\mathrm{d}s}{\mathrm{d}\theta}=|x'(\theta)|=\frac{1}{\kappa(x(\theta))}.
\end{equation}

For each positive solution $v$ of problem \eqref{eq:le} in $\Omega$, define
\begin{equation}\label{eq:k}
k_v(\theta)=\frac{q_v(x(\theta))^2}{2\kappa(x(\theta))}, \quad \theta \in \mathbb{S}^1.
\end{equation}
By \eqref{eq:reg}, we have $k_v>0$ on $\mathbb{S}^1$. Using $\mathrm{d}\theta=\kappa\mathrm{d}s$ in \eqref{eq:flux}, we obtain
\begin{equation*}
\int_0^{2\pi}k_v(\theta)n(\theta)\mathrm{d}\theta=\frac{1}{2}\int_{\partial\Omega}q_v^2\nu\mathrm{d}s=0.
\end{equation*}

We now express the support function of $\overline{\Omega}$ in terms of $x(\theta)$. Since $n(\theta)$ is the outward unit normal at $x(\theta)$, convexity of $\Omega$ gives
\begin{equation*}
(y-x(\theta))\cdot n(\theta)\leq0\quad\text{for every }y\in\overline{\Omega}.
\end{equation*}
Hence the support function of $\overline{\Omega}$ (see \cite[Section 1.7]{schneider}) satisfies
\begin{equation} \label{def-h}
h(\theta)=\max_{y\in\overline{\Omega}}y\cdot n(\theta)=x(\theta)\cdot n(\theta).
\end{equation}
Differentiating with respect to $\theta$ and using \eqref{eq:x'-theta}, we obtain
\begin{equation*}
\begin{aligned}
h'(\theta)&=x'(\theta)\cdot n(\theta)+x(\theta)\cdot\tau(\theta)=x(\theta)\cdot\tau(\theta),\\
h''(\theta)&=x'(\theta)\cdot\tau(\theta)-x(\theta)\cdot n(\theta)=\frac{1}{\kappa(x(\theta))}-h(\theta).
\end{aligned}
\end{equation*}
Thus $h+h''=1/\kappa(x(\theta))$. Since $x\in C^{2,\delta}(\mathbb{S}^1)$ and $h'=x\cdot\tau$, we have $h\in C^{3,\delta}(\mathbb{S}^1)$.

\vskip 0.1cm
Conversely, for $h\in C^{3,\delta}(\mathbb{S}^1)$, set
\begin{equation*}
\Omega_h=\Big\{y\in\mathbb{R}^2:~ y\cdot n(\theta)<h(\theta)\ \text{for every }\theta\in\mathbb{S}^1\Big\}.
\end{equation*}

\begin{lemma}
Assume that $\rho=h+h''>0$ on $\mathbb{S}^1$. Then $\Omega_h$ is the unique bounded convex domain with support function $h$, and $\partial\Omega_h$ is of class $C^{3,\delta}$. For each $\theta\in\mathbb{S}^1$, the boundary point $x_h(\theta)$ with outward unit normal $n(\theta)$ satisfies
\begin{equation}\label{eq:param}
\begin{aligned}
x_h(\theta)&=h(\theta)n(\theta)+h'(\theta)\tau(\theta),\\
x_h'(\theta)&=\rho(\theta)\tau(\theta),\quad \kappa(x_h(\theta))=\frac{1}{\rho(\theta)}.
\end{aligned}
\end{equation}
\end{lemma}

\begin{proof}
Set $x_h=hn+h'\tau$. Since $n'=\tau$ and $\tau'=-n$,
\begin{equation*}
x_h\cdot n=h,\quad x_h'=h'n+h\tau+h''\tau-h'n=\rho\tau.
\end{equation*}
Fix $0<r\leq\min_{\mathbb{S}^1}\rho$ and take $|\theta-\theta_0|\leq\pi$ by periodicity. On the interval $I$ between $\theta_0$ and $\theta$, integrating $x_h'=\rho\tau$ and using $\tau(t)\cdot n(\theta)=\sin(\theta-t)$ give
\begin{equation}\label{eq:support}
\begin{aligned}
h(\theta)-x_h(\theta_0)\cdot n(\theta)
&=\int_{\theta_0}^{\theta}\rho(t)\sin(\theta-t)\mathrm{d}t
=\int_I\rho(t)|\sin(\theta-t)|\mathrm{d}t\\
&\geq r\int_0^{|\theta-\theta_0|}\sin t\mathrm{d}t
=r\left(1-\cos(\theta-\theta_0)\right).
\end{aligned}
\end{equation}
For $|z|<r$, \eqref{eq:support} and $n(\theta_0)\cdot n(\theta)=\cos(\theta-\theta_0)$ imply
\begin{equation*}
\left(x_h(\theta_0)-rn(\theta_0)+z\right)\cdot n(\theta)\leq h(\theta)-r+|z|<h(\theta).
\end{equation*}
By the definition of $\Omega_h$,
\begin{equation}\label{eq:disk}
B_r\left(x_h(\theta_0)-rn(\theta_0)\right)\subset\Omega_h.
\end{equation}
Set $K=\operatorname{conv}(x_h(\mathbb{S}^1))$. For $y\in K$, write
\begin{equation*}
y=\sum_{j=1}^N\lambda_jx_h(\theta_j),\quad \lambda_j>0,\quad \sum_{j=1}^N\lambda_j=1.
\end{equation*}
By \eqref{eq:support}, we have
\begin{equation*}
h(\theta)-y\cdot n(\theta)\geq r\sum_{j=1}^N\lambda_j\left(1-\cos(\theta-\theta_j)\right)\geq0.
\end{equation*}
If $y\cdot n(\theta)=h(\theta)$, then $\theta_j=\theta$ for every $j$, so $y=x_h(\theta)$. Conversely, $x_h(\theta)\cdot n(\theta)=h(\theta)$ by \eqref{eq:param}. By \cite[Section 1.3]{schneider}, we know
\begin{equation*}
\partial K=x_h(\mathbb{S}^1),\quad K=\left\{y\in\mathbb{R}^2 \mid y\cdot n(\theta)\leq h(\theta)\ \text{for every }\theta\in\mathbb{S}^1\right\}.
\end{equation*}
For $y\in K$, equality $y\cdot n(\theta)=h(\theta)$ holds for some $\theta$ if and only if $y\in\partial K$. Thus $\Omega_h=\operatorname{int}K$. By \eqref{eq:disk}, $\Omega_h$ is nonempty. Since $h$ determines $K$, it determines $\Omega_h$ uniquely.

\vskip 0.1cm
Let $s$ denote arc length in the direction of $\tau$. Since $x_h'=\rho\tau$ and $\rho>0$,
\begin{equation*}
\frac{\mathrm{d}s}{\mathrm{d}\theta}=|x_h'(\theta)|=\rho(\theta),\quad \frac{\mathrm{d}}{\mathrm{d}s}x_h(\theta(s))=\tau(\theta(s)).
\end{equation*}
Since $\mathrm{d}s/\mathrm{d}\theta=\rho>0$ and $\rho\in C^{1,\delta}(\mathbb{S}^1)$, the inverse of $\theta\mapsto s(\theta)$ is of class $C^{2,\delta}$. The derivative of the arc length parametrization is $\tau(\theta(s))$, so $\partial\Omega_h$ is of class $C^{3,\delta}$. By \eqref{eq:Fre} and $n'=\tau$,
\begin{equation*}
\kappa(x_h(\theta))=\frac{n'(\theta)\cdot\tau(\theta)}{\mathrm{d}s/\mathrm{d}\theta}=\frac{1}{\rho(\theta)}.
\end{equation*}
\end{proof}

\section{Comparison of positive solutions}\label{sec:comp}

In this section, we compare positive solutions $u$ and $v$ of \eqref{eq:le} in bounded convex domains $\Omega_1,\Omega_2\subset\mathbb{R}^2$ of class $C^{3,\delta}$ with positive boundary curvature. We prove that if $D^2(v^{-\beta})>0$ in $\Omega_2$ and $k_u=k_v$ on $\mathbb{S}^1$, then there exists $c\in\mathbb{R}^2$ such that
\begin{equation*}
\Omega_1=c+\Omega_2,\quad u(x)=v(x-c)\quad\text{for all }x\in\Omega_1.
\end{equation*}
By Lemma \ref{lem:conv}, the condition $D^2(v^{-\beta})>0$ holds when $v$ is the positive ground state of \eqref{eq:le} in $\Omega_2$.

\begin{lemma}\label{lem:harm}
Let $\Omega\subset\mathbb{R}^2$ be a bounded domain of class $C^{3,\delta}$ and fix $2<r<\infty$. Suppose that $U,V\in W^{1,r}(\Omega)$ satisfy
\begin{equation}\label{eq:cr}
\begin{cases}
\partial_{x_1}U=\partial_{x_2}V,\\[2mm]
\partial_{x_2}U=-\partial_{x_1}V
\end{cases}
\quad\text{in }\Omega.
\end{equation}
Assume that $U(x_0)=V(x_0)=0$ for some $x_0\in\Omega$ and that, for some $h\in C^1(\partial\Omega)$,
\begin{equation}\label{eq:uvbd}
U=Ah,\quad V=B\partial_\tau h\quad\text{on }\partial\Omega,
\end{equation}
where $A\in W^{1,r}(\Omega)$ is positive on $\overline\Omega$ and $B\in C(\partial\Omega)$ is nonnegative. Then it holds
$$U=V=0~~~\text{in}~~~ \Omega.$$
\end{lemma}

\begin{proof}
Fix $\zeta\in C^1(\overline\Omega)$ and set $F=(-\partial_{x_2}\zeta,\partial_{x_1}\zeta)$. By \eqref{eq:cr} and the definition of $F$,
\begin{equation*}
\Delta U=\partial_{x_1}\partial_{x_2}V-\partial_{x_2}\partial_{x_1}V=0,\quad \operatorname{div}F=-\partial_{x_1}\partial_{x_2}\zeta+\partial_{x_2}\partial_{x_1}\zeta=0\quad\text{in }\mathcal{D}'(\Omega).
\end{equation*}
Then, $VF$ is continuous on $\overline{\Omega}$ and satisfies
\begin{equation*}
\operatorname{div}(VF)=F\cdot\nabla V\quad\text{in }\mathcal{D}'(\Omega).
\end{equation*}
Since $F\cdot\nu=-\partial_\tau\zeta$ on $\partial\Omega$, the Gauss-Green formula yields
\begin{equation}\label{eq:green}
\int_\Omega\nabla U\cdot\nabla\zeta\mathrm{d}x=\int_\Omega F\cdot\nabla V\mathrm{d}x=-\int_{\partial\Omega}V\partial_\tau\zeta\mathrm{d}s.
\end{equation}

Set $w=\frac{U}{A}$. Since $r>2$, the Sobolev embedding theorem yields $U,A\in C(\overline\Omega)$. In particular, $A_0=\min_{\overline\Omega}A>0$. By the Sobolev chain and product rules and \eqref{eq:uvbd},
\begin{equation}\label{eq:gradw}
w\in W^{1,r}(\Omega),\quad \nabla w=\frac{\nabla U-w\nabla A}{A},\quad w=h\quad\text{on }\partial\Omega.
\end{equation}
For $\varepsilon>0$, define $\eta_\varepsilon(t)=\frac{t}{\sqrt{t^2+\varepsilon^2}}$. Then $|\eta_\varepsilon|\leq1$ and
\begin{equation*}
\eta_\varepsilon'(t)=\frac{\varepsilon^2}{(t^2+\varepsilon^2)^{\frac{3}{2}}}\geq0,\qquad t^2\eta_\varepsilon'(t)\leq\varepsilon.
\end{equation*}
Now using \cite[Chapter 6]{gt}, we can choose $H\in C^1(\overline\Omega)$ with $H=h$ on $\partial\Omega$. Since $w=h$ on $\partial\Omega$, using the trace theorem, we have
\begin{equation*}
\eta_\varepsilon(w)-\eta_\varepsilon(H)\in H^1_0(\Omega).
\end{equation*}
Set $G_\varepsilon=\nabla U\cdot\nabla\eta_\varepsilon(w)$. Testing $\Delta U=0$ with $\eta_\varepsilon(w)-\eta_\varepsilon(H)$, we obtain
\begin{equation*}
\int_\Omega G_\varepsilon\mathrm{d}x=\int_\Omega\nabla U\cdot\nabla\eta_\varepsilon(H)\mathrm{d}x.
\end{equation*}
On $\partial\Omega$, we have
\begin{equation*}
\partial_\tau\eta_\varepsilon(H)=\eta_\varepsilon'(h)\partial_\tau h.
\end{equation*}
Taking $\zeta=\eta_\varepsilon(H)$ in \eqref{eq:green} and using \eqref{eq:uvbd}, we find
\begin{equation}\label{eq:trunc}
\int_\Omega G_\varepsilon\mathrm{d}x=-\int_{\partial\Omega}V\eta_\varepsilon'(h)\partial_\tau h\mathrm{d}s=-\int_{\partial\Omega}B\eta_\varepsilon'(h)|\partial_\tau h|^2\mathrm{d}s\leq0.
\end{equation}

By \eqref{eq:gradw}, we get
\begin{equation*}
\begin{aligned}
G_\varepsilon&=\frac{\eta_\varepsilon'(w)}{A}\left(|\nabla U|^2-w\nabla U\cdot\nabla A\right) =\frac{\eta_\varepsilon'(w)}{A}\left(\left|\nabla U-\frac{w}{2}\nabla A\right|^2-\frac{w^2}{4}|\nabla A|^2\right).
\end{aligned}
\end{equation*}
Since $w^2\eta_\varepsilon'(w)\leq\varepsilon$ and $A\geq A_0$, it follows that
\begin{equation*}
G_\varepsilon\geq-\frac{\varepsilon}{4A_0}|\nabla A|^2,\qquad |G_\varepsilon|\leq G_\varepsilon+\frac{\varepsilon}{2A_0}|\nabla A|^2.
\end{equation*}
Integration and \eqref{eq:trunc} imply
\begin{equation}\label{eq:grad}
\|G_\varepsilon\|_{L^1(\Omega)}\leq\int_\Omega G_\varepsilon\mathrm{d}x+\frac{\varepsilon}{2A_0}\int_\Omega|\nabla A|^2\mathrm{d}x\leq\frac{\varepsilon}{2A_0}\int_\Omega|\nabla A|^2\mathrm{d}x.
\end{equation}

Fix $\varphi\in C_c^\infty(\Omega)$. Testing $\Delta U=0$ with $\varphi\eta_\varepsilon(w)$ and using \eqref{eq:grad}, we obtain
\begin{equation*}
\left|\int_\Omega\eta_\varepsilon(w)\nabla U\cdot\nabla\varphi\mathrm{d}x\right|
=\left|\int_\Omega\varphi G_\varepsilon\mathrm{d}x\right|\leq\|\varphi\|_{L^\infty(\Omega)}\|G_\varepsilon\|_{L^1(\Omega)}\to0.
\end{equation*}
Since $w=\frac{U}{A}$ and $A>0$, we have $\eta_\varepsilon(w)\to\operatorname{sgn}(U)$ a.e. in $\Omega$. As $|\eta_\varepsilon|\leq1$ and $|\nabla U|\in L^1(\Omega)$, the dominated convergence theorem and the chain rule give
\begin{equation*}
\eta_\varepsilon(w)\nabla U\to\operatorname{sgn}(U)\nabla U=\nabla|U|\quad\text{in }L^1(\Omega;\mathbb{R}^2).
\end{equation*}
Hence
\begin{equation*}
\int_\Omega\nabla|U|\cdot\nabla\varphi\mathrm{d}x=0\quad\text{for every }\varphi\in C_c^\infty(\Omega).
\end{equation*}
By Weyl's lemma, $|U|$ is harmonic in $\Omega$. Since $|U|$ attains its minimum $0$ at $x_0\in\Omega$, the strong minimum principle yields $U\equiv0$. Then $\nabla V=0$ by \eqref{eq:cr}, so $V\equiv V(x_0)=0$ in $\Omega$.
\end{proof}

For a matrix field $T=(T_{ij})$, we write
\begin{equation*}
(\operatorname{div}T)_i=\partial_{x_1}T_{i1}+\partial_{x_2}T_{i2}\quad i=1,2,\qquad |T|^2=\sum_{i,j=1}^2T_{ij}^2.
\end{equation*}

\begin{lemma}\label{lem:tensor}
Let $\Omega\subset\mathbb{R}^2$ be a bounded convex domain of class $C^{3,\delta}$. Let $T\in C(\overline\Omega;\mathbb{R}^{2\times2})$ be symmetric, with $\operatorname{div}T\in L^\infty(\Omega;\mathbb{R}^2)$. Assume that, for some $C>0$,
\begin{equation}\label{eq:divt}
\operatorname{tr}T=0\quad\text{in }\Omega,\qquad |\operatorname{div}T|\leq C|T|\quad\text{a.e. in }\Omega.
\end{equation}
Suppose that
\begin{equation}\label{eq:tbd}
T(\nu,\nu)=ah,\quad T(\nu,\tau)=-b\partial_\tau h\quad\text{on }\partial\Omega,
\end{equation}
where $a,h\in C^1(\partial\Omega)$, $b\in C(\partial\Omega)$ and $a,b>0$ on $\partial\Omega$. Then $T=0$ in $\Omega$.
\end{lemma}

\begin{proof}
The equations \eqref{eq:cr} can be written as
\begin{equation*}
\partial_{x_1}\begin{pmatrix}U\\[1mm]V\end{pmatrix}
=\begin{pmatrix}\partial_{x_2}V\\[1mm]-\partial_{x_2}U\end{pmatrix}
=-J\partial_{x_2}\begin{pmatrix}U\\[1mm]V\end{pmatrix},
\quad J=\begin{pmatrix}0&-1\\1&0\end{pmatrix}.
\end{equation*}
Since $T$ is symmetric and $\operatorname{tr}T=0$, we have $T_{21}=T_{12}$ and $T_{22}=-T_{11}$. Set $Q=(T_{11},-T_{12})^{\mathsf{T}}$. Then $|T|^2=2|Q|^2$ and
\begin{equation*}
(\partial_{x_1}+J\partial_{x_2})Q
=\begin{pmatrix}\partial_{x_1}T_{11}+\partial_{x_2}T_{12}\\[1mm]-\partial_{x_1}T_{12}+\partial_{x_2}T_{11}\end{pmatrix}
=\begin{pmatrix}(\operatorname{div}T)_1\\[1mm]-(\operatorname{div}T)_2\end{pmatrix}
\quad\text{in }\mathcal{D}'(\Omega).
\end{equation*}

After a translation, we assume that $0\in\Omega$. Denote
\begin{equation*}
Z=P(x)Q, \quad P(x)=(x_1^2-x_2^2)I+2x_1x_2J.
\end{equation*}
Using $J^2=-I$ and $J^{\mathsf T}=-J$, we have
\begin{equation*}
PJ=JP,\quad \partial_{x_1}P+J\partial_{x_2}P=0,\quad P^{\mathsf T}P=|x|^4I.
\end{equation*}
Thus $Z(0)=0$, and $|T|^2=2|Q|^2$ yields
\begin{equation*}
2|Z|^2=2|x|^4|Q|^2=|x|^4|T|^2.
\end{equation*}
Set $F=(\partial_{x_1}+J\partial_{x_2})Z$. Using $JP=PJ$ and the product rule in $\mathcal{D}'(\Omega)$, we obtain
\begin{equation*}
F=(\partial_{x_1}P+J\partial_{x_2}P)Q+P(\partial_{x_1}+J\partial_{x_2})Q=P\begin{pmatrix}(\operatorname{div}T)_1\\[1mm]-(\operatorname{div}T)_2\end{pmatrix}.
\end{equation*}
Since $|P(x)\xi|=|x|^2|\xi|$ for every $\xi\in\mathbb{R}^2$, it follows from \eqref{eq:divt} that
\begin{equation}\label{eq:fbound}
|F|=|x|^2|\operatorname{div}T|\leq C|x|^2|T|=\sqrt{2}C|Z|\quad\text{a.e. in }\Omega.
\end{equation}
As $Z\in C(\overline\Omega;\mathbb{R}^2)$ and $\Omega$ is bounded, we have $F\in L^\infty(\Omega;\mathbb{R}^2)$.

\vspace{0.1cm}

Applying $\partial_{x_1}-J\partial_{x_2}$ to $F=(\partial_{x_1}+J\partial_{x_2})Z$, we obtain in $\mathcal{D}'(\Omega)$
\begin{equation*}
\begin{aligned}
\partial_{x_1}F-\partial_{x_2}(JF)
&=\partial_{x_1x_1}Z+J\left(\partial_{x_1x_2}Z-\partial_{x_2x_1}Z\right)-J^2\partial_{x_2x_2}Z=\partial_{x_1x_1}Z+\partial_{x_2x_2}Z=\Delta Z.
\end{aligned}
\end{equation*}
Fix $2<r<\infty$. Since $F,JF\in L^r(\Omega;\mathbb{R}^2)$, we have
\begin{equation*}
\partial_{x_1}F-\partial_{x_2}(JF)\in W^{-1,r}(\Omega;\mathbb{R}^2).
\end{equation*}
Using \cite[Chapter 4]{simader}, we obtain a unique $Y\in W^{1,r}_0(\Omega;\mathbb{R}^2)$ satisfying
\begin{equation*}
\begin{cases}
\Delta Y=\partial_{x_1}F-\partial_{x_2}(JF) &\text{in }\Omega,\\[1mm]
Y=0 &\text{on }\partial\Omega.
\end{cases}
\end{equation*}
Since $\Delta(Z-Y)=0$ in $\mathcal{D}'(\Omega)$, Weyl's lemma yields $Z-Y\in C^\infty(\Omega;\mathbb{R}^2)$. Thus
\begin{equation*}
Z\in W^{1,r}_{\mathrm{loc}}(\Omega;\mathbb{R}^2).
\end{equation*}

On $\{Z\neq0\}$, the vectors $Z$ and $JZ$ are orthogonal and have the same length. Define
\begin{equation*}
c_1=\frac{F\cdot Z}{|Z|^2},\quad c_2=\frac{F\cdot JZ}{|Z|^2}.
\end{equation*}
Then $F=c_1Z+c_2JZ$ almost everywhere on $\{Z\neq0\}$, and \eqref{eq:fbound} implies
\begin{equation*}
c_1^2+c_2^2=\frac{|F|^2}{|Z|^2}\leq 2C^2.
\end{equation*}
Set $c_1=c_2=0$ on $\{Z=0\}$, where $F=0$ a.e. by \eqref{eq:fbound}. Then $c_1,c_2\in L^\infty(\Omega)$ and
\begin{equation}\label{eq:sys}
(\partial_{x_1}+J\partial_{x_2})Z=(c_1I+c_2J)Z\quad\text{a.e. in }\Omega.
\end{equation}

Let $d=\operatorname{dist}(0,\partial\Omega)>0$. For each $x\in\partial\Omega$, convexity of $\Omega$ gives $(y-x)\cdot\nu(x)\leq0$ for all $y\in\overline\Omega$. Taking $y=d\nu(x)\in\overline{B_d(0)}\subset\overline\Omega$, we obtain $x\cdot\nu(x)\geq d$. Thus
\begin{equation*}
\chi(x)=\arctan\left(\frac{Jx\cdot\nu(x)}{x\cdot\nu(x)}\right)
\in C^1(\partial\Omega).
\end{equation*}
Set $R(t)=(\cos t)I+(\sin t)J$. In the orthonormal basis $\frac{x}{|x|},\frac{Jx}{|x|}$, the components of $\nu(x)$ are
\begin{equation*}
\cos\chi(x)=\frac{x\cdot\nu(x)}{|x|},\qquad \sin\chi(x)=\frac{Jx\cdot\nu(x)}{|x|}.
\end{equation*}
Hence $R(\chi(x))x=|x|\nu(x)$. Since $P(R(t)x)=R(2t)P(x)$, we have
\begin{equation*}
R(2\chi(x))P(x)=|x|^2P(\nu(x))\quad\text{on }\partial\Omega.
\end{equation*}
By \cite[Chapter 6]{gt}, we can extend $-2\chi$ to a function $\psi_0\in C^1(\overline\Omega)$. Since $c_1,c_2\in L^\infty(\Omega)$, we have
\begin{equation*}
(c_2,-c_1)-\nabla\psi_0\in L^r(\Omega;\mathbb{R}^2).
\end{equation*}
By \cite[Chapter 4]{simader}, there exists $z\in W^{1,r}_0(\Omega)$ satisfying
\begin{equation*}
\Delta z=\operatorname{div}\left((c_2,-c_1)-\nabla\psi_0\right)\quad\text{in }\mathcal{D}'(\Omega).
\end{equation*}
Then $\psi=z+\psi_0\in W^{1,r}(\Omega)$ satisfies
\begin{equation}\label{eq:psi}
\begin{cases}
\Delta\psi=\partial_{x_1}c_2-\partial_{x_2}c_1 &\text{in }\mathcal{D}'(\Omega),\\
\psi=-2\chi &\text{on }\partial\Omega.
\end{cases}
\end{equation}

The vector field $(c_1+\partial_{x_2}\psi,c_2-\partial_{x_1}\psi)$ belongs to $L^r(\Omega;\mathbb{R}^2)$. By \eqref{eq:psi},
\begin{equation*}
\partial_{x_1}(c_2-\partial_{x_1}\psi)-\partial_{x_2}(c_1+\partial_{x_2}\psi)
=\partial_{x_1}c_2-\partial_{x_2}c_1-\Delta\psi=0
\quad\text{in }\mathcal{D}'(\Omega).
\end{equation*}
Since $\Omega$ is convex and hence simply connected, then by \cite{mardare}, we have $\varphi\in W^{1,r}(\Omega)$ such that
\begin{equation}\label{eq:phi}
\partial_{x_1}\varphi=c_1+\partial_{x_2}\psi,\qquad \partial_{x_2}\varphi=c_2-\partial_{x_1}\psi.
\end{equation}
Set $S=\mathrm{e}^{-\varphi}R(-\psi)$ and $(U,V)^{\mathsf T}=SZ$. Since $\varphi,\psi\in W^{1,r}(\Omega)$ and $r>2$, the Sobolev embedding theorem gives $\varphi,\psi\in C(\overline\Omega)$. By the Sobolev chain rule, $S\in W^{1,r}(\Omega;\mathbb{R}^{2\times2})$. Since $Z\in C(\overline\Omega;\mathbb{R}^2)\cap W^{1,r}_{\mathrm{loc}}(\Omega;\mathbb{R}^2)$, we obtain
\begin{equation*}
U,V\in C(\overline\Omega)\cap W^{1,r}_{\mathrm{loc}}(\Omega).
\end{equation*}
Using $R'=JR$ and $RJ=JR$, we have
\begin{equation*}
\partial_{x_j}S=-S\left((\partial_{x_j}\varphi)I+(\partial_{x_j}\psi)J\right),\quad j=1,2.
\end{equation*}
Thus $J^2=-I$ and \eqref{eq:phi} yield
\begin{equation*}
\begin{aligned}
\partial_{x_1}S+J\partial_{x_2}S
&=-S\left((\partial_{x_1}\varphi-\partial_{x_2}\psi)I+(\partial_{x_1}\psi+\partial_{x_2}\varphi)J\right) =-S(c_1I+c_2J).
\end{aligned}
\end{equation*}
Since $SJ=JS$, the product rule and \eqref{eq:sys} give
\begin{equation*}
\begin{aligned}
(\partial_{x_1}+J\partial_{x_2})(SZ)
&=(\partial_{x_1}S+J\partial_{x_2}S)Z+S(\partial_{x_1}+J\partial_{x_2})Z\\
&=-S(c_1I+c_2J)Z+S(c_1I+c_2J)Z=0.
\end{aligned}
\end{equation*}
Writing this equation in components and using $Z(0)=0$, we obtain
\begin{equation}\label{eq:uv}
\begin{cases}
\partial_{x_1}U=\partial_{x_2}V,\\
\partial_{x_2}U=-\partial_{x_1}V
\end{cases}
\quad\text{a.e. in }\Omega,\qquad U(0)=V(0)=0.
\end{equation}

On $\partial\Omega$, we have $\psi=-2\chi$ and $R(2\chi)P(x)=|x|^2P(\nu)$. Since $Z=P(x)Q$, it follows that
\begin{equation*}
\begin{pmatrix}U\\V\end{pmatrix}=|x|^2\mathrm{e}^{-\varphi}P(\nu)Q\quad\text{on }\partial\Omega.
\end{equation*}
Using $Q=(T_{11},-T_{12})^{\mathsf T}$, $T_{22}=-T_{11}$ and $\tau=(-\nu_2,\nu_1)$, we find
\begin{equation*}
\begin{aligned}
(P(\nu)Q)_1&=(\nu_1^2-\nu_2^2)T_{11}+2\nu_1\nu_2T_{12}=T(\nu,\nu),
\end{aligned}
\end{equation*}
and\begin{equation*}
\begin{aligned}
(P(\nu)Q)_2&=2\nu_1\nu_2T_{11}-(\nu_1^2-\nu_2^2)T_{12}=-T(\nu,\tau).
\end{aligned}
\end{equation*}
The boundary conditions \eqref{eq:tbd} therefore yield
\begin{equation}\label{eq:uvtr}
U=|x|^2\mathrm{e}^{-\varphi}ah,\quad V=|x|^2\mathrm{e}^{-\varphi}b\partial_\tau h\quad\text{on }\partial\Omega.
\end{equation}
Since $|x|\geq d>0$ and $a>0$ on $\partial\Omega$, we have $\log(|x|^2a)\in C^1(\partial\Omega)$. By \cite[Chapter 6]{gt}, we choose $g,H\in C^1(\overline\Omega)$ satisfying
\begin{equation*}
g=\log(|x|^2a),\quad H=h\quad\text{on }\partial\Omega.
\end{equation*}
Set $A=\mathrm{e}^{g-\varphi}$ on $\overline\Omega$ and $B=|x|^2\mathrm{e}^{-\varphi}b$ on $\partial\Omega$. Then $A\in W^{1,r}(\Omega)$ is positive on $\overline\Omega$, and $B\in C(\partial\Omega)$ is positive. By \eqref{eq:uvtr}, we have
\begin{equation*}
U=Ah,\quad V=B\partial_\tau h\quad\text{on }\partial\Omega.
\end{equation*}

Since $AH\in W^{1,r}(\Omega)$, we have $-\Delta(AH)\in W^{-1,r}(\Omega)$. By \cite[Chapter 4]{simader}, there exists $w\in W^{1,r}_0(\Omega)$ satisfying
\begin{equation*}
\Delta w=-\Delta(AH)\quad\text{in }\mathcal{D}'(\Omega).
\end{equation*}
Then $U_0=AH+w$ belongs to $W^{1,r}(\Omega)\cap C(\overline\Omega)$ and satisfies
\begin{equation*}
\begin{cases}
\Delta U_0=0 &\text{in }\mathcal{D}'(\Omega),\\
U_0=Ah &\text{on }\partial\Omega.
\end{cases}
\end{equation*}
By \eqref{eq:uv}, we have
$$\Delta U=\partial_{x_1x_2}V-\partial_{x_2x_1}V=0~~~in ~~~\mathcal{D}'(\Omega).$$ Hence $U-U_0$ is harmonic by Weyl's lemma and vanishes on $\partial\Omega$. The maximum principle yields $U=U_0$, so $U\in W^{1,r}(\Omega)$. Using \eqref{eq:uv} again, we obtain
\begin{equation*}
\nabla V=(-\partial_{x_2}U,\partial_{x_1}U)\in L^r(\Omega;\mathbb{R}^2).
\end{equation*}
Since $V\in C(\overline\Omega)$, we also have $V\in W^{1,r}(\Omega)$. Lemma \ref{lem:harm} gives $U=V=0$ in $\Omega$. Thus $SZ=0$, and the invertibility of $S$ yields $Z=0$. Since $|x|^4|T(x)|^2=2|Z(x)|^2$, we have $T(x)=0$ for $x\neq0$. By continuity, $T(0)=0$.
\end{proof}

Let $v$ be a positive solution of problem \eqref{eq:le} in $\Omega_2$. Set $\alpha=\frac{2}{p-1}$. For $s>0$ and $c\in\mathbb{R}^2$, define
\begin{equation*}
w_{s,c}(z)=s^{-\alpha}v\left(\frac{z-c}{s}\right),\quad z\in c+s\Omega_2.
\end{equation*}
Since $\alpha p=\alpha+2$, we have
\begin{equation*}
-\Delta w_{s,c}(z)=s^{-\alpha-2}v\left(\frac{z-c}{s}\right)^p=w_{s,c}(z)^p.
\end{equation*}
Thus $w_{s,c}$ is a positive solution of problem \eqref{eq:le} in $c+s\Omega_2$.

\begin{lemma}\label{lem:match}
Assume that $D^2(v^{-\frac{p-1}{2}})>0$ in $\Omega_2$. For each $x\in\mathbb{R}^2$, $a>0$ and $b\in\mathbb{R}^2$, the equations
\begin{equation*}
w_{s,c}(x)=a,\quad \nabla w_{s,c}(x)=b
\end{equation*}
have a unique solution $(s,c)\in(0,\infty)\times\mathbb{R}^2$ with $x\in c+s\Omega_2$. The Hessian $M(a,b):=D^2w_{s,c}(x)$ depends only on $(a,b)$, and
\begin{equation*}
M\in C^\infty\big((0,\infty)\times\mathbb{R}^2;\mathbb{R}^{2\times2}\big),\quad \operatorname{tr}M(a,b)=-a^p.
\end{equation*}
\end{lemma}

\begin{proof}
Set $G=v^{-\frac{p-1}{2}}$. For distinct $y_1,y_2\in\Omega_2$, set $\xi=y_1-y_2$. Since $\Omega_2$ is convex and $D^2G>0$ in $\Omega_2$,
\begin{equation*}
\left(\nabla G(y_1)-\nabla G(y_2)\right)\cdot\xi=\int_0^1D^2G(y_2+t\xi)[\xi,\xi]\mathrm{d}t>0.
\end{equation*}
Thus $\nabla G$ is injective. Fix $z\in\mathbb{R}^2$. Since $\Omega_2$ is bounded and $v=0$ on $\partial\Omega_2$,
\begin{equation*}
G(y)-z\cdot y\to+\infty\quad\text{as }y\to\partial\Omega_2.
\end{equation*}
Hence $G(y)-z\cdot y$ attains its minimum at some $y_z\in\Omega_2$, where $\nabla G(y_z)=z$. Thus $\nabla G$ is both injective and surjective. Since $v>0$ in $\Omega_2$, using the Schauder estimates on $-\Delta v=v^p$, we get $v\in C^\infty(\Omega_2)$. As $D^2G$ is positive definite, the inverse function theorem shows that $\nabla G:\Omega_2\to\mathbb{R}^2$ is a smooth diffeomorphism.

\vskip 0.1cm

Fix $x\in\mathbb{R}^2$, $a>0$ and $b\in\mathbb{R}^2$. At $x=c+sy$ with $y\in\Omega_2$, the identity $\frac{(p-1)\alpha}{2} =1$ gives
\begin{equation*}
w_{s,c}(x)^{-\frac{p-1}{2}}=sG(y),\quad \nabla(w_{s,c}^{-\frac{p-1}{2}})(x)=\nabla G(y).
\end{equation*}
By the chain rule, the conditions $w_{s,c}(x)=a$ and $\nabla w_{s,c}(x)=b$ are equivalent to
\begin{equation*}
\nabla G(y)=-\frac{p-1}{2} a^{-\frac{p+1}{2}}b,\quad
s=\frac{a^{-\frac{p-1}{2}}}{G(y)}.
\end{equation*}
The first equation has a unique solution $y\in\Omega_2$ because $\nabla G$ is bijective. For this $y$, the second equation determines $s>0$, and $x=c+sy$ determines $c$. Differentiating $w_{s,c}$ twice, we obtain
\begin{equation}\label{eq:match}
\begin{aligned}
y&=(\nabla G)^{-1}(-\frac{p-1}{2} a^{-\frac{p+1}{2}}b),\quad s=\frac{a^{-\frac{p-1}{2}}}{G(y)}, \quad
M(a,b) =D^2w_{s,c}(x)=s^{-\alpha-2}D^2v(y).
\end{aligned}
\end{equation}
It follows from \eqref{eq:match} that $M$ is smooth on $(0,\infty)\times\mathbb{R}^2$. At $x=c+sy$, the definition of $w_{s,c}$ and the equation $-\Delta w_{s,c}=w_{s,c}^p$ give us that
\begin{equation}\label{eq:mid}
\begin{aligned}
a&=s^{-\alpha}v(y),\quad b=s^{-\alpha-1}\nabla v(y), \quad
\operatorname{tr}M(a,b) =\Delta w_{s,c}(x)=-a^p.
\end{aligned}
\end{equation}
For fixed $s>0$ and $c\in\mathbb{R}^2$, applying \eqref{eq:match} at each $x\in c+s\Omega_2$ with $a=w_{s,c}(x)$ and $b=\nabla w_{s,c}(x)$, we have
\begin{equation}\label{eq:m}
D^2w_{s,c}=M(w_{s,c},\nabla w_{s,c})\quad\text{in }c+s\Omega_2.
\end{equation}
\end{proof}

Assume that $D^2(v^{-\frac{p-1}{2}})>0$ in $\Omega_2$. Let $u$ be a positive solution of problem \eqref{eq:le} in $\Omega_1$. For each $x\in\Omega_1$, choose $(s,c)$ as in Lemma \ref{lem:match} with $a=u(x)$ and $b=\nabla u(x)$. As in \cite{mira}, set
\begin{equation}\label{eq:t}
T(x)=D^2u(x)-D^2w_{s,c}(x)=D^2u(x)-M(u(x),\nabla u(x)).
\end{equation}
For $i=1,2$, let $\kappa_i$ denote the curvature of $\partial\Omega_i$. The point $x_i(\theta)$ is the unique point of $\partial\Omega_i$ with outward unit normal $n(\theta)$. With $q_1=q_u$ and $q_2=q_v$, we write
\begin{equation*}
q_i(\theta)=q_i(x_i(\theta)),\quad \kappa_i(\theta)=\kappa_i(x_i(\theta)),\quad k_i(\theta)=\frac{q_i(\theta)^2}{2\kappa_i(\theta)}.
\end{equation*}
Primes denote differentiation with respect to $\theta$. Set $\gamma=\frac{1}{\alpha+1}$ and $m=2-\gamma$.

\begin{lemma}\label{lem:field}
The map $M$ defined in Lemma \ref{lem:match} has a $C^1$ extension to
\begin{equation*}
\Big\{(a,b)\in\mathbb{R}\times\mathbb{R}^2:a>0\ \text{or}\ b\neq0\Big\}.
\end{equation*}
Moreover, $T\in C^1(\overline\Omega_1;\mathbb{R}^{2\times2})$, $\operatorname{tr}T=0$, and
\begin{equation*}
(\operatorname{div}T)_i=-\sum_{j,\ell=1}^2\partial_{b_\ell}M_{ij}(u,\nabla u)T_{\ell j}\quad\text{in }\Omega_1,\quad i=1,2.
\end{equation*}
In particular, $|\operatorname{div}T|\leq C|T|$ in $\Omega_1$. At $x_1(\theta)$, we have
\begin{equation}\label{eq:bd}
T(n,n)=\frac{q_1^{1+\gamma}}{2}\left(\frac{q_1^m}{k_1}-\frac{q_2^m}{k_2}\right),\quad
T(n,\tau)=-\frac{q_1^{1+\gamma}}{2m}\left(\frac{(q_1^m)'}{k_1}-\frac{(q_2^m)'}{k_2}\right).
\end{equation}
\end{lemma}

\begin{proof}
By \eqref{eq:reg} and \cite[Chapter 6]{gt}, we can extend $v$ as a $C^{3,\delta}$ function to a neighborhood of $\overline\Omega_2$. For $s>0$ and $y$ in this neighborhood, set
\begin{equation*}
\Psi(s,y)=\left(s^{-\alpha}v(y),s^{-\alpha-1}\nabla v(y)\right).
\end{equation*}
Fix $b_0\in\mathbb{R}^2\setminus\{0\}$. We solve $\Psi(s_0,y_0)=(0,b_0)$ for $s_0>0$ and $y_0\in\partial\Omega_2$. Since $v=0$ and $\nabla v=-q_2\nu$ on $\partial\Omega_2$, this equation is equivalent to
\begin{equation*}
b_0=-s_0^{-\alpha-1}q_2(y_0)\nu(y_0).
\end{equation*}
As $q_2(y_0)>0$, the direction of $b_0$ determines $\nu(y_0)=-\frac{b_0}{|b_0|}$. The outward normal map is bijective, so $y_0$ is uniquely determined. Using $\gamma=\frac{1}{\alpha+1}$, we obtain
\begin{equation*}
|b_0|=s_0^{-\alpha-1}q_2(y_0),\quad s_0=\left(\frac{q_2(y_0)}{|b_0|}\right)^\gamma.
\end{equation*}

Let $(\sigma,\eta)\in\ker D\Psi(s_0,y_0)$ and write $\nu=\nu(y_0)$ and $\tau=\tau(y_0)$. Using $v(y_0)=0$ and $\nabla v(y_0)=-q_2(y_0)\nu$ in $D\Psi(s_0,y_0)(\sigma,\eta)=0$, we obtain
\begin{equation}\label{eq:kerpsi}
\begin{cases}
\nu\cdot\eta=0,\\[1mm]
D^2v(y_0)\eta+\frac{\alpha+1}{s_0}q_2(y_0)\sigma\nu=0.
\end{cases}
\end{equation}
By the first equation in \eqref{eq:kerpsi}, we get
\begin{equation*}
\eta=(\eta\cdot\tau)\tau.
\end{equation*}
Since $\nu\cdot\tau=0$, from \eqref{eq:hess} and the second equation in \eqref{eq:kerpsi}, we can deduce
\begin{equation*}
0=(D^2v(y_0)\eta)\cdot\tau=(\eta\cdot\tau)D^2v(y_0)(\tau,\tau)=-\kappa_2(y_0)q_2(y_0)(\eta\cdot\tau).
\end{equation*}
Since $\kappa_2(y_0)q_2(y_0)>0$, we have $\eta=(\eta\cdot\tau)\tau=0$. The second equation in \eqref{eq:kerpsi} then gives $\sigma=0$, so $D\Psi(s_0,y_0)$ is invertible.

\vspace{0.1cm}

By the inverse function theorem, the equation $\Psi(s,y)=(a,b)$ defines $C^1$ functions $s=s(a,b)>0$ and $y=y(a,b)$ near $(0,b_0)$, with
\begin{equation*}
s(0,b_0)=s_0,\quad y(0,b_0)=y_0.
\end{equation*}
For $z\notin\overline{\Omega}_2$ sufficiently close to $y_0$, let $\xi\in\partial\Omega_2$ be the nearest point to $z$. By \cite[Lemma 14.16]{gt},
\begin{equation*}
z=\xi+t\nu(\xi),\quad t=|z-\xi|=\operatorname{dist}(z,\partial\Omega_2).
\end{equation*}
As $z\to y_0$, we have $\xi\to y_0$ and $t\to0$. Since $v(\xi)=0$ and $\partial_\nu v(y_0)=-q_2(y_0)<0$, by the continuity of $\nabla v$ and $\nu$, we have
\begin{equation*}
v(z)=\int_0^t\nabla v(\xi+r\nu(\xi))\cdot\nu(\xi)\mathrm{d}r\leq-\frac{q_2(y_0)}{2}t<0.
\end{equation*}
For $a>0$, the first component of $\Psi(s,y)=(a,b)$ gives $v(y)=s^\alpha a>0$. Thus $y\in\Omega_2$.
At $x=sy\in s\Omega_2$, the definition of $w_{s,0}$ implies
\begin{equation*}
w_{s,0}(x)=s^{-\alpha}v(y)=a,\quad \nabla w_{s,0}(x)=s^{-\alpha-1}\nabla v(y)=b.
\end{equation*}
Hence, by Lemma \ref{lem:match}, we have
\begin{equation*}
M(a,b)=D^2w_{s,0}(x)=s(a,b)^{-\alpha-2}D^2v(y(a,b)).
\end{equation*}
The last expression is $C^1$ near $(0,b_0)$ since $s_0>0$ and $v\in C^{3,\delta}$ near $y_0$. It therefore extends $M$ across $(0,b_0)$, with
\begin{equation}\label{eq:mbd}
M(0,b_0)=s_0^{-\alpha-2}D^2v(y_0).
\end{equation}
Fix $b\neq0$. Since $M$ has a local $C^1$ extension near $(0,b)$, the limits of $M(a,b)$ and its first derivatives exist as $a\to0^+$. Set
\begin{equation*}
M(0,b)=\lim_{t\to0^+}M(t,b),\quad M(a,b)=2M(0,b)-M(-a,b)\quad\text{for }a<0.
\end{equation*}
For $a<0$ and $b\neq0$,
\begin{equation*}
\partial_a M(a,b)=\partial_a M(-a,b),\quad D_bM(a,b)=2D_bM(0,b)-D_bM(-a,b).
\end{equation*}
Thus $M$ is $C^1$ across $a=0$ for every $b\neq0$.
By \eqref{eq:reg}, $u\in C^{3,\delta}(\overline\Omega_1)$ and $\nabla u=-q_1\nu\neq0$ on $\partial\Omega_1$. Since $u>0$ in $\Omega_1$, $M$ is $C^1$ near $(u(x),\nabla u(x))$ for every $x\in\overline\Omega_1$. It follows from \eqref{eq:t} that
\begin{equation*}
T\in C^1(\overline\Omega_1;\mathbb{R}^{2\times2}),\quad \sup_{\overline\Omega_1}|D_bM(u,\nabla u)|<\infty.
\end{equation*}
By \eqref{eq:mid}, $\operatorname{tr}T=\Delta u+u^p=0$. Fix $x_0\in\Omega_1$. Choose $w=w_{s,c}$ by Lemma \ref{lem:match} so that
\begin{equation*}
w(x_0)=u(x_0)=a,\quad \nabla w(x_0)=\nabla u(x_0)=b.
\end{equation*}
With $s$ and $c$ fixed, \eqref{eq:m} gives $D^2w=M(w,\nabla w)$ near $x_0$. Hence it holds
\begin{equation*}
\begin{aligned}
D^2(u-w)(x_0)&=D^2u(x_0)-M(a,b)=T(x_0),
\end{aligned}
\end{equation*}
and
\begin{equation*}
\begin{aligned}
\nabla\Delta(u-w)(x_0)&=\nabla(w^p-u^p)(x_0)=-pa^{p-1}\nabla(u-w)(x_0)=0.
\end{aligned}
\end{equation*}
By \eqref{eq:t} and the chain rule, for $i=1,2$, we get
\begin{equation}\label{eq:div}
\begin{aligned}
(\operatorname{div}T)_i(x_0)
&=\partial_{x_i}\Delta(u-w)(x_0)
-\sum_{j=1}^2\partial_{x_j}\left(M_{ij}(u,\nabla u)-M_{ij}(w,\nabla w)\right)(x_0)\\
&=-\sum_{j=1}^2\partial_aM_{ij}(a,b)\partial_{x_j}(u-w)(x_0)-\sum_{j,\ell=1}^2\partial_{b_\ell}M_{ij}(a,b)\partial_{x_\ell x_j}(u-w)(x_0)\\
&=-\sum_{j,\ell=1}^2\partial_{b_\ell}M_{ij}(a,b)T_{\ell j}(x_0).
\end{aligned}
\end{equation}
Since $D_bM(u,\nabla u)$ is bounded on $\overline\Omega_1$, \eqref{eq:div} gives  $|\operatorname{div}T|\leq C|T|$ in $\Omega_1$.

\vskip 0.1cm

Fix $\theta$ and write $x_i=x_i(\theta)$, $n=n(\theta)$ and $\tau=\tau(\theta)$. Choose $s>0$ so that $s^{-\alpha-1}q_2=q_1$. Recall $\gamma=\frac{1}{\alpha+1}$, it follows that
\begin{equation}\label{eq:scale}
s=\left(\frac{q_2}{q_1}\right)^\gamma,\quad s^{-\alpha-2}=\left(\frac{q_1}{q_2}\right)^{1+\gamma}.
\end{equation}
Using $v(x_2)=0$ and $\nabla v(x_2)=-q_2n$,
\begin{equation*}
\Psi(s,x_2)=(0,-s^{-\alpha-1}q_2n)=(0,-q_1n)=(u(x_1),\nabla u(x_1)).
\end{equation*}
By \eqref{eq:t} and \eqref{eq:mbd}, we have
\begin{equation*}
T(x_1)=D^2u(x_1)-M(0,-q_1n)=D^2u(x_1)-s^{-\alpha-2}D^2v(x_2).
\end{equation*}
Since $\mathrm{d}\theta=\kappa_i\mathrm{d}s_i$, we have $(\partial_\tau q_i)(x_i)=\kappa_iq_i'$. Applying \eqref{eq:hess} to $u$ at $x_1$ and to $v$ at $x_2$, we obtain
\begin{equation*}
T(n,n)=\kappa_1q_1-s^{-\alpha-2}\kappa_2q_2,\quad T(n,\tau)=-\kappa_1q_1'+s^{-\alpha-2}\kappa_2q_2'.
\end{equation*}
Substituting \eqref{eq:scale} and $\kappa_i=\frac{q_i^2}{2k_i}$, and using $m=2-\gamma$ and $(q_i^m)'=mq_i^{m-1}q_i'$, we find
\begin{equation*}
T(n,n)=\frac{q_1^{1+\gamma}}{2}\left(\frac{q_1^m}{k_1}-\frac{q_2^m}{k_2}\right),\quad T(n,\tau)=-\frac{q_1^{1+\gamma}}{2m}\left(\frac{(q_1^m)'}{k_1}-\frac{(q_2^m)'}{k_2}\right).
\end{equation*}
\end{proof}

We now use the boundary identities \eqref{eq:bd} to prove the rigidity result.

\begin{proposition}\label{prop:comp}
Let $\Omega_1,\Omega_2\subset\mathbb{R}^2$ be bounded convex domains of class $C^{3,\delta}$ with positive boundary curvature. Suppose that $u$ and $v$ are positive solutions of problem \eqref{eq:le} in $\Omega_1$ and $\Omega_2$, respectively, and that
\begin{equation}\label{eq:hyp}
D^2(v^{-\frac{p-1}{2}})>0\quad\text{in }\Omega_2,\quad k_u=k_v\quad\text{on }\mathbb{S}^1.
\end{equation}
Then there exists $c\in\mathbb{R}^2$ such that
\begin{equation*}
\Omega_1=c+\Omega_2,\quad u(x)=v(x-c)\quad\text{for all }x\in\Omega_1.
\end{equation*}
\end{proposition}

\begin{proof}
By \eqref{eq:hyp}, $k_1=k_2$. Denote their common value by $k$. Then
\begin{equation}\label{eq:equal}
k(\theta)=\frac{q_1(\theta)^2}{2\kappa_1(\theta)}=\frac{q_2(\theta)^2}{2\kappa_2(\theta)}.
\end{equation}
On $\partial\Omega_1$, set
$H(x_1(\theta))=q_1(\theta)^m-q_2(\theta)^m$.
Since $H\in C^1(\partial\Omega_1)$ and $\mathrm{d}\theta=\kappa_1\mathrm{d}s_1$, we have
\begin{equation*}
\partial_\tau H(x_1(\theta))=\kappa_1(\theta)\frac{\mathrm{d}}{\mathrm{d}\theta}\left(q_1(\theta)^m-q_2(\theta)^m\right).
\end{equation*}
Since $2k\kappa_1=q_1^2$ by \eqref{eq:equal}, the boundary identities \eqref{eq:bd} become
\begin{equation}\label{eq:tnn}
\begin{aligned}
T(\nu,\nu)&=\frac{q_1^{1+\gamma}}{2k}H=aH,\quad a=\kappa_1q_1^{\gamma-1},\\
T(\nu,\tau)&=-\frac{q_1^{1+\gamma}}{2mk\kappa_1}\partial_\tau H=-b\partial_\tau H,\quad b=\frac{q_1^{\gamma-1}}{m}.
\end{aligned}
\end{equation}
Here, $a\in C^1(\partial\Omega_1)$ and $b\in C(\partial\Omega_1)$ are positive. By \eqref{eq:t}, $T$ is symmetric. Lemma \ref{lem:field} yields $T\in C^1(\overline\Omega_1;\mathbb{R}^{2\times2})$ and
\begin{equation*}
\operatorname{tr}T=0,\quad |\operatorname{div}T|\leq C|T|\quad\text{in }\Omega_1.
\end{equation*}
Applying Lemma \ref{lem:tensor} with $h=H$, we obtain $T=0$ in $\Omega_1$. Since $T$ is continuous on $\overline\Omega_1$,
\begin{equation}\label{eq:tzero}
T=0\quad\text{on }\overline\Omega_1.
\end{equation}
By \eqref{eq:tnn} and \eqref{eq:tzero}, we get $aH=0$ on $\partial\Omega_1$. Since $a>0$, we have $H=0$. Thus $q_1^m=q_2^m$, and hence $q_1=q_2$ on $\mathbb{S}^1$. It follows from \eqref{eq:equal} that
\begin{equation*}
\kappa_1(\theta)=\frac{q_1(\theta)^2}{2k(\theta)}=\frac{q_2(\theta)^2}{2k(\theta)}=\kappa_2(\theta).
\end{equation*}
By \eqref{eq:x'-theta},
\begin{equation*}
(x_1-x_2)'(\theta)=\left(\frac{1}{\kappa_1(\theta)}-\frac{1}{\kappa_2(\theta)}\right)\tau(\theta)=0.
\end{equation*}
Hence $x_1(\theta)=c+x_2(\theta)$ for some fixed $c\in\mathbb{R}^2$ and every $\theta$, so $\partial\Omega_1=c+\partial\Omega_2$. By \cite[Section 1.4]{schneider},
\begin{equation*}
\overline{\Omega_1}=\operatorname{conv}(\partial\Omega_1)=c+\operatorname{conv}(\partial\Omega_2)=c+\overline{\Omega_2}.
\end{equation*}
Taking interiors implies $\Omega_1=c+\Omega_2$.

\vskip 0.1cm

Set $\widetilde{v}(x)=v(x-c)$ and $w=u-\widetilde{v}$ in $\Omega_1$. Both $u$ and $\widetilde{v}$ vanish on $\partial\Omega_1$. Since $x_1(\theta)-c=x_2(\theta)$ and the outward normal to $\partial\Omega_i$ at $x_i(\theta)$ is $n(\theta)$,
\begin{equation*}
\partial_\nu w(x_1(\theta))=\left(\nabla u(x_1(\theta))-\nabla v(x_2(\theta))\right)\cdot n(\theta)=-q_1(\theta)+q_2(\theta)=0.
\end{equation*}
Subtracting the equations for $u$ and $\widetilde{v}$, we obtain
\begin{equation*}
\begin{cases}
-\Delta w=A(x)w &\text{in }\Omega_1,\\
w=\partial_\nu w=0 &\text{on }\partial\Omega_1,
\end{cases}
\end{equation*}
where
\begin{equation*}
A(x)=p\int_0^1\big((1-t)\widetilde{v}(x)+tu(x)\big)^{p-1}\mathrm{d}t.
\end{equation*}
Since $u$ and $\widetilde{v}$ are bounded, $A\in L^\infty(\Omega_1)$. By the unique continuation principle, we deduce that $w=0$ in $\Omega_1$. Thus $u(x)=v(x-c)$ for every $x\in\Omega_1$.
\end{proof}

\section{Prescribing the boundary function}\label{sec:real}

In this section, we fix a function $k\in C^2(\mathbb{S}^1)$ satisfying
\begin{equation}\label{eq:target}
k>0\quad\text{on }\mathbb{S}^1,\quad \int_0^{2\pi}k(\theta)n(\theta)\mathrm{d}\theta=0.
\end{equation}
Fix $0<\delta<\min\{1,p-1\}$. We construct a bounded convex domain $\Omega$ with $C^{3,\delta}$ boundary of positive curvature. The positive ground state $v$ of \eqref{eq:le} in $\Omega$ satisfies $k_v=k$ on $\mathbb{S}^1$. The domain $\Omega$ is unique up to translation.

\vskip 0.1cm

Set
\begin{equation}\label{eq:def-x+y}
\begin{aligned}
X&=\left\{h\in C^{3,\delta}(\mathbb{S}^1) :~\int_0^{2\pi}h(\theta)n(\theta)\mathrm{d}\theta=0\right\}, \\
Y&=\left\{k\in C^{1,\delta}(\mathbb{S}^1) :~ \int_0^{2\pi}k(\theta)n(\theta)\mathrm{d}\theta=0\right\},
\end{aligned}
\end{equation}
and
\begin{equation}\label{eq:def-ma}
\mathcal{A}=\big\{h\in X :~ h+h''>0\text{ on }\mathbb{S}^1\big\}.
\end{equation}
Let $h\in\mathcal{A}$. Since $h\in X$ and $\tau'=-n$, integration by parts in \eqref{eq:param} implies
\begin{equation*}
\int_0^{2\pi}x_h(\theta_0)\mathrm{d}\theta_0
=\int_0^{2\pi}(hn+h'\tau)\mathrm{d}\theta_0=\int_0^{2\pi}(hn-h\tau')\mathrm{d}\theta_0
=2\int_0^{2\pi}hn\mathrm{d}\theta_0=0.
\end{equation*}
Fix $0<r\leq\min_{\mathbb{S}^1}(h+h'')$. For each $\theta\in\mathbb{S}^1$, integrating \eqref{eq:support} with respect to $\theta_0$ yields
\begin{equation*}
2\pi h(\theta)=\int_0^{2\pi}\left(h(\theta)-x_h(\theta_0)\cdot n(\theta)\right)\mathrm{d}\theta_0\geq r\int_0^{2\pi}\left(1-\cos(\theta-\theta_0)\right)\mathrm{d}\theta_0=2\pi r.
\end{equation*}
Thus $h\geq r>0$. If $|y|<r$, then
\begin{equation*}
y\cdot n(\theta)\leq |y|<r\leq h(\theta) \quad \text{for every } \theta\in\mathbb{S}^1.
\end{equation*}
Hence $B_r(0)\subset\Omega_h$ by the definition of $\Omega_h$.

\vskip 0.1cm
For $h\in C^{3,\delta}(\mathbb{S}^1)$ with $h+h''>0$, let $v_h$ be the positive ground state of problem \eqref{eq:le} in $\Omega_h$. Define
\begin{equation}\label{eq:map}
q_h(\theta)=-\nabla v_h(x_h(\theta))\cdot n(\theta),\quad \mathcal{K}(h)=\frac{1}{2}q_h^2(h+h'').
\end{equation}
By \eqref{eq:reg} and \eqref{eq:param}, $q_h\in C^{2,\delta}(\mathbb{S}^1)$. Since $\kappa(x_h(\theta))=1/(h(\theta)+h''(\theta))$, \eqref{eq:k} and \eqref{eq:map} give $\mathcal{K}(h)=k_{v_h}\in C^{1,\delta}(\mathbb{S}^1)$. Using $\mathrm{d}s=(h+h'')\mathrm{d}\theta$ in \eqref{eq:flux}, we obtain
\begin{equation*}
\int_0^{2\pi}\mathcal{K}(h)(\theta)n(\theta)\mathrm{d}\theta=\frac{1}{2}\int_{\partial\Omega_h}(\partial_\nu v_h)^2\nu\mathrm{d}s=0.
\end{equation*}
Thus $\mathcal{K}(h)\in Y$.

\begin{lemma}\label{lem:pull}
Let $\Omega\subset\mathbb{R}^2$ be a bounded domain and $\Phi$ a $C^2$ diffeomorphism defined on a neighborhood of $\overline{\Omega}$. Write $\Omega_\Phi=\Phi(\Omega)$ and $A=(D\Phi)^{-1}$.
\begin{enumerate}
\item[(i)] The Sobolev constant satisfies
\begin{equation}\label{eq:scont}
S_p(\Omega_\Phi)\to S_p(\Omega)\quad\text{as }\Phi\to\mathrm{Id}\text{ in }C^1(\overline{\Omega}).
\end{equation}

\item[(ii)] For $1<r<\infty$ and $z\in W^{2,r}(\Omega)$,
\begin{equation}\label{eq:pull}
\big[\Delta(z\circ\Phi^{-1})\big]\circ\Phi=\sum_{i,j,\ell=1}^2 A_{ji}\partial_{x_j}\left(A_{\ell i}\partial_{x_\ell}z\right)\quad\text{a.e. in }\Omega.
\end{equation}
For $V\in C_c^2(\mathbb{R}^2;\mathbb{R}^2)$ and $|t|$ sufficiently small, set
\begin{equation*}
\Phi_t=\mathrm{Id}+tV,\quad \mathcal{L}_t z=[-\Delta(z\circ\Phi_t^{-1})]\circ\Phi_t.
\end{equation*}
If $v$ is a positive solution of problem \eqref{eq:le} in $\Omega$, then
\begin{equation}\label{eq:opder}
\mathcal{L}'_0v=-L_v(\nabla v\cdot V)\quad\text{in }\Omega.
\end{equation}
\end{enumerate}
\end{lemma}

\begin{proof}
For $\Phi$ sufficiently close to $\mathrm{Id}$ in $C^1(\overline{\Omega})$, we have $J_\Phi=\det D\Phi>0$. For $w\in H^1_0(\Omega)$, the change of variables $y=\Phi(x)$ gives
\begin{equation}\label{eq:change}
\int_{\Omega_\Phi}|\nabla(w\circ\Phi^{-1})|^2\mathrm{d}y=\int_\Omega|A^{\mathsf{T}}\nabla w|^2J_\Phi\mathrm{d}x,\quad
\int_{\Omega_\Phi}|w\circ\Phi^{-1}|^{p+1}\mathrm{d}y=\int_\Omega|w|^{p+1}J_\Phi\mathrm{d}x.
\end{equation}
As $\Phi\to\mathrm{Id}$ in $C^1(\overline{\Omega})$, then
\begin{equation*}
J_\Phi=1+o(1),\quad |A^{\mathsf{T}}\xi|^2J_\Phi=1+o(1)
\end{equation*}
uniformly for $x\in\overline{\Omega}$ and $|\xi|=1$. By \eqref{eq:change}, we get
\begin{equation*}
\frac{\|\nabla(w\circ\Phi^{-1})\|_{L^2(\Omega_\Phi)}^2}{\|w\circ\Phi^{-1}\|_{L^{p+1}(\Omega_\Phi)}^2}=(1+o(1))\frac{\|\nabla w\|_{L^2(\Omega)}^2}{\|w\|_{L^{p+1}(\Omega)}^2}
\end{equation*}
uniformly for $0\ne w\in H^1_0(\Omega)$. Since $w\mapsto w\circ\Phi^{-1}$ is a bijection between $H^1_0(\Omega)$ and $H^1_0(\Omega_\Phi)$, taking the infimum proves \eqref{eq:scont}.

For $z\in W^{2,r}(\Omega)$, write $y=\Phi(x)$. Since $D\Phi^{-1}(y)=A(x)$, the chain rule gives
\begin{equation*}
\big[\partial_{y_i}(z\circ\Phi^{-1})\big]\circ\Phi=\sum_{\ell=1}^2 A_{\ell i}\partial_{x_\ell}z.
\end{equation*}
Applying the same chain rule to $\partial_{y_i}(z\circ\Phi^{-1})$, we obtain
\begin{equation*}
[\partial_{y_i y_i}(z\circ\Phi^{-1})]\circ\Phi
=\sum_{j=1}^2 A_{ji}\partial_{x_j} \big([\partial_{y_i}(z\circ\Phi^{-1})]\circ\Phi \big)=\sum_{j,\ell=1}^2 A_{ji}\partial_{x_j}\left(A_{\ell i}\partial_{x_\ell}z\right),
\end{equation*}
which yields \eqref{eq:pull} by summing.

\vskip 0.1cm
Fix an open set $U\subset\subset\Omega$. For $|t|$ sufficiently small, $\Phi_t(\overline{U})\subset\Omega$. Set $v_t=v\circ\Phi_t$ on $U$. By the definition of $\mathcal{L}_t$,
\begin{equation*}
\mathcal{L}_t v_t=(-\Delta v)\circ\Phi_t=v_t^p\quad\text{in }U.
\end{equation*}
Since $\Phi_t(x)=x+tV(x)$ and $v\in C^3(\Omega)$, we have
\begin{equation*}
\frac{v_t(x)-v(x)}{t} = \int_0^1\nabla v(x+\sigma tV(x))\cdot V(x)\mathrm{d}\sigma ~~ \to ~~ \nabla v(x)\cdot V(x)\quad\text{in }C^2(\overline{U}).
\end{equation*}
Since $(D\Phi_t)^{-1}$ depends smoothly on $t$ in $C^1(\overline{U})$, \eqref{eq:pull} shows that $t\mapsto\mathcal{L}_t$ is differentiable in the operator norm from $C^2(\overline{U})$ to $C(\overline{U})$. In particular, $\mathcal{L}_t\to\mathcal{L}_0=-\Delta$ in such norms. Hence,
\begin{equation*}
\begin{aligned}
pv^{p-1}(\nabla v\cdot V)
& = \lim_{t\to0}\frac{\mathcal{L}_t v_t - \mathcal{L}_0 v}{t}\\
& = \lim_{t\to0} \frac{\mathcal{L}_t v - \mathcal{L}_0 v}{t}
+\lim_{t\to0}\mathcal{L}_t \Big(\frac{v_t - v}{t}\Big)\\
&=\mathcal{L}'_0v-\Delta(\nabla v\cdot V).
\end{aligned}
\end{equation*}
Thus $\mathcal{L}'_0v=-L_v(\nabla v\cdot V)$ in $U\subset\subset\Omega$, which proves \eqref{eq:opder}.
\end{proof}

Set $F(t)=(t_+)^p$ for $t\in\mathbb{R}$, where $t_+=\max\{t,0\}$. Then $F\in C^1(\mathbb{R})$ and $F'(t)=p(t_+)^{p-1}$.

\begin{lemma}\label{lem:quot}
Let $\Omega\subset\mathbb{R}^2$ be a bounded domain of class $C^{2,\delta}$. Set
\begin{equation*}
E=C^{2,\delta}_0(\overline{\Omega})=\left\{z\in C^{2,\delta}(\overline{\Omega}):~ z=0\ \text{on }\partial\Omega\right\}
\end{equation*}
and equip $E$ with the $C^{2,\delta}(\overline{\Omega})$ norm. Assume that $v\in E$, $v>0$ in $\Omega$ and $\partial_\nu v<0$ on $\partial\Omega$. Then $z/v$ extends to a function in $C^{1,\delta}(\overline{\Omega})$ for every $z\in E$, with
\begin{equation}\label{ee:est-z/v}
\left\|\frac{z}{v}\right\|_{C^{1,\delta}(\overline{\Omega})}\leq C\|z\|_{C^{2,\delta}(\overline{\Omega})}.
\end{equation}
Moreover, the map $z\mapsto F(z)$ is $C^\infty$ from a neighborhood of $v$ in $E$ to $C^{0,\delta}(\overline{\Omega})$.
\end{lemma}

\begin{proof}
Write $d(x)=\operatorname{dist}(x,\partial\Omega)$ for $x\in\overline{\Omega}$. For $s_0>0$ sufficiently small and $0<d(x)<s_0$, let $\xi$ be the nearest point on $\partial\Omega$ to $x$. By \cite[Lemma 14.16]{gt}, $x=\xi-s\nu(\xi)$ with $s=d(x)$.
Since $z(\xi)=0$,
\begin{equation}\label{eq:z/d=}
\frac{z(\xi-s\nu(\xi))}{s}=-\int_0^1\nabla z(\xi-ts\nu(\xi))\cdot\nu(\xi)\mathrm{d}t.
\end{equation}
At $s=0$, the right-hand side of \eqref{eq:z/d=} equals $-\partial_\nu z(\xi)$.
Since $\nabla z\in C^{1,\delta}(\overline{\Omega})$ and the normal coordinates are of class $C^{1,\delta}$, the chain rule yields
\begin{equation*}
\left.\frac{z}{d}\right|_{\partial\Omega}=-\partial_\nu z,\quad
\left\|\frac{z}{d}\right\|_{C^{1,\delta}(\{d\leq s_0/2\})}\leq C\|\nabla z\|_{C^{1,\delta}(\overline{\Omega})}.
\end{equation*}
Taking $z=v$, we have $(v/d)|_{\partial\Omega}=-\partial_\nu v>0$.
Then, by continuity of $v/d$ and compactness of $\partial\Omega$, there exists $0<d_0<s_0/2$ such that
\begin{equation*}
\frac{v}{d}\geq\frac{1}{2} \min_{\partial\Omega}(-\partial_\nu v)>0, \quad \text{on }\{d\leq d_0\}.
\end{equation*}
On $\{d\leq d_0\}$, the quotient $(z/d)/(v/d)$ extends $z/v$ to the boundary with value $\partial_\nu z/\partial_\nu v$. On $\{d\geq d_0/2\}$, $v$ has a positive minimum.
Thus, it follows \eqref{ee:est-z/v}.

\vskip 0.1cm
By \eqref{ee:est-z/v}, the map $\eta\mapsto\eta/v$ is bounded and linear from $E$ to $C^{1,\delta}(\overline{\Omega})$. For $\|\eta\|_E$ sufficiently small, $1+\eta/v\geq1/2$ and
\begin{equation*}
F(v+\eta)=v^p\left(1+\frac{\eta}{v}\right)^p
=v^p\sum_{j=0}^{\infty}\binom{p}{j}\left(\frac{\eta}{v}\right)^j.
\end{equation*}
The series for $(1+\eta/v)^p$ converges in $C^{1,\delta}(\overline{\Omega})$. Multiplication by $v^p\in C^{0,\delta}(\overline{\Omega})$ gives a convergent power series for $F(v+\eta)$ in $C^{0,\delta}(\overline{\Omega})$.
\end{proof}

\begin{lemma}\label{lem:shape}
The following statements hold.
\begin{enumerate}
\item[(i)] The maps
\begin{equation*}
h\mapsto q_h\in C^{1,\delta}(\mathbb{S}^1),\quad h\mapsto\mathcal{K}(h)\in Y
\end{equation*}
are $C^\infty$ on the open set of $h\in C^{3,\delta}(\mathbb{S}^1)$ satisfying $h+h''>0$.

\item[(ii)] Fix $h$ as in (i) and set
\begin{equation*}
v=v_h,\quad q=q_h,\quad \rho=h+h'',\quad \kappa=\frac{1}{\rho}.
\end{equation*}
For $f\in C^{3,\delta}(\mathbb{S}^1)$, let $w\in C^{2,\delta}(\overline{\Omega_h})$ be the unique solution of
\begin{equation}\label{eq:w}
\begin{cases}
L_vw=0 &\text{in }\Omega_h,\\
w(x_h(\theta))=q(\theta)f(\theta) &\text{for }\theta\in\mathbb{S}^1.
\end{cases}
\end{equation}
Writing $Dq_h[f]$ for the derivative of $h\mapsto q_h$ in the direction $f$, we have
\begin{equation}\label{eq:dq}
Dq_h[f]=-(\partial_\nu w)\circ x_h-\kappa qf+\kappa q'f'
\end{equation}
and
\begin{equation}\label{eq:dk}
D\mathcal{K}(h)[f]=\frac{q^2}{2}(f+f'')+q\rho Dq_h[f].
\end{equation}
\end{enumerate}
\end{lemma}

\begin{proof}
Fix $h$ as in (i) and set $\Omega=\Omega_h$. By \eqref{eq:reg} and \eqref{eq:param}, $v\in C^{3,\delta}(\overline{\Omega})$ and $q\in C^{2,\delta}(\mathbb{S}^1)$. Take $0<s_0<\min_{\mathbb{S}^1}\rho$ sufficiently small. For $x\in\mathbb{R}^2$ with $\operatorname{dist}(x,\partial\Omega)<s_0$, let $x_h(\theta)$ be the nearest point on $\partial\Omega$ to $x$. By \cite[Lemma 14.16]{gt},
\begin{equation*}
x=x_h(\theta)+sn(\theta),\quad |s|=\operatorname{dist}(x,\partial\Omega),
\end{equation*}
where $s<0$ in $\Omega$ and $s>0$ outside $\overline{\Omega}$. The map $(\theta,s)\mapsto x_h(\theta)+sn(\theta)$ is of class $C^{2,\delta}$. By \eqref{eq:param} and $\det(\tau,n)=-1$,
\begin{equation*}
\det\left(\partial_\theta(x_h+sn),\partial_s(x_h+sn)\right)=\det((\rho+s)\tau,n)=-(\rho+s)\neq0.
\end{equation*}
Thus the inverse map $x\mapsto(\theta,s)$ is of class $C^{2,\delta}$.

\vskip 0.1cm
Choose $\chi\in C_c^\infty((-s_0,s_0))$ with $\chi=1$ near $0$.
For $f\in C^{3,\delta}(\mathbb{S}^1)$, we set
\begin{equation*}
V_f(x_h(\theta)+sn(\theta))=\chi(s)\left(f(\theta)n(\theta)+f'(\theta)\tau(\theta)\right),\quad |s|<s_0,
\end{equation*}
and extend $V_f$ by zero outside this neighborhood. Then, the map $f\mapsto V_f$ is linear and satisfies the estimate
\begin{equation*}
\|V_f\|_{C^{2,\delta}(\mathbb{R}^2)}\leq C\|f\|_{C^{3,\delta}(\mathbb{S}^1)}.
\end{equation*}
Set $\Phi_f=\mathrm{Id}+V_f$.
For $\|f\|_{C^{3,\delta}(\mathbb{S}^1)}$ sufficiently small, $\rho+f+f''>0$ and $\|DV_f\|_{L^\infty}<1/2$. For each $y\in\mathbb{R}^2$, the equation $x=y-V_f(x)$ has a unique solution by the contraction mapping theorem.
Thus $\Phi_f$ is bijective.
Since $I+DV_f$ is invertible at every point, $\Phi_f$ is a $C^{2,\delta}$ diffeomorphism of $\mathbb{R}^2$. By $\chi(0)=1$ and \eqref{eq:param}, we get
\begin{equation*}
\Phi_f(x_h(\theta)) = x_h(\theta)+ f(\theta)n(\theta)+ f'(\theta)\tau(\theta)=x_{h+f}(\theta).
\end{equation*}
Hence $\Phi_f(\partial\Omega)=\partial\Omega_{h+f}$. Since $\Phi_f(\Omega)$ is the bounded component of $\mathbb{R}^2\setminus\partial\Omega_{h+f}$, we have $\Phi_f(\Omega)=\Omega_{h+f}$.

\vskip 0.1cm
Write $E=C^{2,\delta}_0(\overline{\Omega})$ and set
\begin{equation*}
A_f=(D\Phi_f)^{-1},\quad \mathcal{L}_fz=[-\Delta(z\circ\Phi_f^{-1})]\circ\Phi_f.
\end{equation*}
Since $f\mapsto DV_f$ is bounded and linear from $C^{3,\delta}(\mathbb{S}^1)$ to $C^{1,\delta}(\overline{\Omega})$, for $f$ sufficiently close to zero, we have
\begin{equation*}
A_f=(I+DV_f)^{-1}=\sum_{j=0}^{\infty}(-DV_f)^j\quad\text{in }C^{1,\delta}(\overline{\Omega}).
\end{equation*}
Thus $A_f$ depends smoothly on $f$. By \eqref{eq:pull}, $\mathcal{L}_f$ also depends smoothly on $f$ in the operator norm from $E$ to $C^{0,\delta}(\overline{\Omega})$. Since $v>0$ in $\Omega$ and $\partial_\nu v<0$ on $\partial\Omega$ by \eqref{eq:reg}, Lemma \ref{lem:quot} shows that
\begin{equation*}
(f,z)\mapsto\mathcal{L}_fz-F(z)
\end{equation*}
is smooth near $(0,v)$ with values in $C^{0,\delta}(\overline{\Omega})$.

\vskip 0.1cm
By \cite[Chapter 6]{gt}, $-\Delta:E\to C^{0,\delta}(\overline{\Omega})$ is an isomorphism. Since $pv^{p-1}\in C^{0,\delta}(\overline{\Omega})$ and $E\hookrightarrow C^{0,\delta}(\overline{\Omega})$ is compact, multiplication by $pv^{p-1}$ is compact from $E$ to $C^{0,\delta}(\overline{\Omega})$. Thus $L_v=-\Delta-pv^{p-1}$ is Fredholm of index zero. Since $E\subset H^2(\Omega)\cap H^1_0(\Omega)$, \eqref{eq:spec} gives $\ker L_v=\{0\}$ in $E$. Then $L_v$ is an isomorphism from $E$ to $C^{0,\delta}(\overline{\Omega})$. At $(f,z)=(0,v)$, we have
\begin{equation*}
\mathcal{L}_0v-F(v)=0,\quad D_z(\mathcal{L}_fz-F(z))\big|_{(0,v)}[\zeta]=L_v\zeta\quad\text{for }\zeta\in E.
\end{equation*}
By the implicit function theorem, for each $f$ near zero, there exists a unique solution $z(f)\in E$ near $v$ of $\mathcal{L}_fz=F(z)$. The map $f\mapsto z(f)$ is smooth.

\vskip 0.1cm
Set
\begin{equation*}
z_f=v_{h+f}\circ\Phi_f,\quad J_f=\det D\Phi_f.
\end{equation*}
By \eqref{eq:reg} and the definition of $\mathcal{L}_f$, we have $z_f\in E$ and $\mathcal{L}_fz_f=F(z_f)$. Since $v_{h+f}$ is a ground state of \eqref{eq:le} in $\Omega_{h+f}$, the change of variables in \eqref{eq:change} and the identities \eqref{eq:energy} and \eqref{eq:least} give
\begin{equation*}
\int_\Omega |A_f^{\mathsf{T}}\nabla z_f|^2J_f\mathrm{d}x=\int_\Omega z_f^{p+1}J_f\mathrm{d}x=S_p(\Omega_{h+f})^{\frac{p+1}{p-1}}.
\end{equation*}
As $f\to0$ in $C^{3,\delta}(\mathbb{S}^1)$, $A_f\to I$ and $J_f\to1$ uniformly on $\overline{\Omega}$. It follows from \eqref{eq:scont} that
\begin{equation*}
\|z_f\|_{H^1_0(\Omega)}\leq C,\quad \int_\Omega z_f^{p+1}\mathrm{d}x\to S_p(\Omega)^{\frac{p+1}{p-1}}.
\end{equation*}

Fix $2<r<\infty$ and set $R_f=\mathcal{L}_f+\Delta$. By \eqref{eq:pull} and $A_f\to I$ in $C^{1,\delta}(\overline{\Omega})$, $R_f\to0$ in the operator norm from $W^{2,r}(\Omega)$ to $L^r(\Omega)$ and from $E$ to $C^{0,\delta}(\overline{\Omega})$. The equation for $z_f$ is
\begin{equation}\label{ee:zf=}
\begin{cases}
-\Delta z_f=F(z_f)-R_fz_f &\text{in }\Omega,\\
z_f=0 &\text{on }\partial\Omega.
\end{cases}
\end{equation}
The Sobolev embedding theorem in dimension two and the Dirichlet $W^{2,r}$ estimate \cite[Chapter 9]{gt} give
\begin{equation*}
\begin{aligned}
\|F(z_f)\|_{L^r(\Omega)}&=\|z_f\|_{L^{pr}(\Omega)}^p\leq C\|z_f\|_{H^1_0(\Omega)}^p\leq C,\\
\|z_f\|_{W^{2,r}(\Omega)}&\leq C\|F(z_f)\|_{L^r(\Omega)}+o(1)\|z_f\|_{W^{2,r}(\Omega)}.
\end{aligned}
\end{equation*}
Absorbing the last term gives a uniform bound in $W^{2,r}(\Omega)$. Since $r>2$, every sequence $f_j\to0$ has a subsequence along which
\begin{equation*}
z_{f_j}\rightharpoonup z\quad\text{in }W^{2,r}(\Omega),\quad z_{f_j}\to z\quad\text{in }C^1(\overline{\Omega}).
\end{equation*}
In \eqref{ee:zf=}, $R_{f_j}z_{f_j}\to0$ in $L^r(\Omega)$ and $F(z_{f_j})\to z^p$ uniformly. Therefore,
\begin{equation*}
\begin{cases}
-\Delta z=z^p,\quad z\geq0 &\text{in }\Omega,\\
z=0 &\text{on }\partial\Omega.
\end{cases}
\end{equation*}
Since $z^p\in C^{0,\delta}(\overline{\Omega})$, the Schauder estimate \cite[Chapter 6]{gt} yields $z\in E$. Passing to the limit in the integral of $z_{f_j}^{p+1}$, we obtain
\begin{equation*}
\int_\Omega z^{p+1}\mathrm{d}x=S_p(\Omega)^{\frac{p+1}{p-1}}>0.
\end{equation*}
The strong maximum principle gives $z>0$ in $\Omega$. By \eqref{eq:least}, $z$ is a ground state of \eqref{eq:le} in $\Omega$, so $z=v$ by Proposition \ref{prop:bf}. Thus $z_f\to v$ in $C^1(\overline{\Omega})$.

\vskip 0.1cm
Since $F\in C^1(\mathbb{R})$, the chain rule gives $F(z_f)\to F(v)$ in $C^1(\overline{\Omega})$ and hence in $C^{0,\delta}(\overline{\Omega})$. Subtracting the equations for $z_f$ and $v$, we obtain
\begin{equation*}
\begin{cases}
\mathcal{L}_f(z_f-v)=F(z_f)-F(v)-R_fv &\text{in }\Omega,\\
z_f-v=0 &\text{on }\partial\Omega.
\end{cases}
\end{equation*}
By \eqref{eq:pull}, the ellipticity bound and the $C^{0,\delta}$ bounds for the coefficients of $\mathcal{L}_f$ are uniform for small $f$. The global Schauder estimate \cite[Chapter 6]{gt} yields
\begin{equation*}
\begin{aligned}
\|z_f-v\|_E
&\leq C\|z_f-v\|_{C(\overline{\Omega})}+C\|F(z_f)-F(v)\|_{C^{0,\delta}(\overline{\Omega})} +C\|R_fv\|_{C^{0,\delta}(\overline{\Omega})}\to0.
\end{aligned}
\end{equation*}
Consequently,
\begin{equation}\label{eq:vcont}
\|v_{h+f}\circ\Phi_f-v_h\|_{C^{2,\delta}(\overline{\Omega_h})}\to0\quad\text{as }f\to0\text{ in }C^{3,\delta}(\mathbb{S}^1).
\end{equation}
Since $\mathcal{L}_fz_f=F(z_f)$ and $z_f\to v$ in $E$ by \eqref{eq:vcont}, local uniqueness gives $z_f=z(f)$ for small $f$. Hence $f\mapsto z_f$ is smooth in $E$.

\vskip 0.1cm
The outward unit normal at $x_{h+f}(\theta)=\Phi_f(x_h(\theta))$ is $n(\theta)$. Differentiating $z_f=v_{h+f}\circ\Phi_f$ with respect to $x$, we obtain
\begin{equation}\label{eq:qpull}
q_{h+f}(\theta)=-n(\theta)\cdot A_f(x_h(\theta))^{\mathsf{T}}\nabla z_f(x_h(\theta)).
\end{equation}
Both $A_f$ and $\nabla z_f$ depend smoothly on $f$ in $C^{1,\delta}(\overline{\Omega})$. Since $x_h$ is fixed, \eqref{eq:qpull} yields smooth dependence of $q_{h+f}$ on $f$ in $C^{1,\delta}(\mathbb{S}^1)$. By \eqref{eq:map},
\begin{equation*}
\mathcal{K}(h+f)=\frac{1}{2}q_{h+f}^2(\rho+f+f'').
\end{equation*}
It follows that $f\mapsto\mathcal{K}(h+f)$ is smooth with values in the closed subspace $Y$.

\vskip 0.1cm
Fix $f\in C^{3,\delta}(\mathbb{S}^1)$ and write $z_t=v_{h+tf}\circ\Phi_{tf}$. Set
\begin{equation*}
\zeta_f=\left.\frac{\mathrm{d}z_t}{\mathrm{d}t}\right|_{t=0},\quad w=\zeta_f-\nabla v\cdot V_f.
\end{equation*}
Since $\zeta_f\in E$ and $v\in C^{3,\delta}(\overline{\Omega})$, we have $w\in C^{2,\delta}(\overline{\Omega})$. Differentiating $\mathcal{L}_{tf}z_t=F(z_t)$ at $t=0$ and applying \eqref{eq:opder} with $\Phi_{tf}=\mathrm{Id}+tV_f$, we obtain
\begin{equation*}
L_v\zeta_f=-\left.\frac{\mathrm{d}}{\mathrm{d}t}(\mathcal{L}_{tf}v)\right|_{t=0}=L_v(\nabla v\cdot V_f).
\end{equation*}
Hence $L_vw=0$ in $\Omega$. Since $z_t=0$ on $\partial\Omega$ for every $t$, $\zeta_f=0$ on $\partial\Omega$. At $x_h(\theta)$,
\begin{equation*}
w=-\nabla v\cdot V_f=qn\cdot(fn+f'\tau)=qf.
\end{equation*}
This proves \eqref{eq:w}. Uniqueness follows from $\ker L_v=\{0\}$ in $E$.

\vskip 0.1cm
Since $\zeta_f=w+\nabla v\cdot V_f$,
\begin{equation*}
\nabla\zeta_f=\nabla w+D^2v V_f+(DV_f)^{\mathsf{T}}\nabla v.
\end{equation*}
Differentiating $A_{tf}(I+tDV_f)=I$ at $t=0$ gives $\left.\frac{\mathrm{d}}{\mathrm{d}t}A_{tf}\right|_{t=0}=-DV_f$. By \eqref{eq:qpull} and $V_f=fn+f'\tau$ at $x_h(\theta)$,
\begin{equation}\label{eq:qdot-hess}
\begin{aligned}
(Dq_h[f])(\theta)&=-n\cdot\left(\nabla\zeta_f-(DV_f)^{\mathsf{T}}\nabla v\right)\\
&=-\partial_\nu w-D^2v(n,V_f)=-\partial_\nu w-fD^2v(n,n)-f'D^2v(n,\tau).
\end{aligned}
\end{equation}
Since $q=q_v\circ x_h$ and $x_h'=\rho\tau$, we have $q'(\theta)=\rho(\theta)\partial_\tau q_v(x_h(\theta))$. Thus \eqref{eq:hess} and $\kappa=1/\rho$ yield at $x_h(\theta)$
\begin{equation*}
D^2v(n,n)=\kappa q,\quad D^2v(n,\tau)=-\partial_\tau q_v=-\kappa q'.
\end{equation*}
Substituting these identities into \eqref{eq:qdot-hess} gives \eqref{eq:dq}. Finally, differentiating
\begin{equation*}
\mathcal{K}(h+tf)=\frac{1}{2}q_{h+tf}^2\left(\rho+t(f+f'')\right)
\end{equation*}
at $t=0$ yields \eqref{eq:dk}.
\end{proof}

\begin{lemma}\label{lem:inv}
For $h\in\mathcal{A}$, the kernel of $D\mathcal{K}(h)$ on $C^{3,\delta}(\mathbb{S}^1)$ is $\operatorname{span}\{\cos\theta,\sin\theta\}$. The restriction of $D\mathcal{K}(h)$ to $X$ is an isomorphism from $X$ to $Y$.
\end{lemma}
\begin{proof}
Fix $h\in\mathcal A$. Write
\begin{equation*}
\Omega=\Omega_h,\quad v=v_h,\quad q=q_h,\quad \rho=h+h'',\quad \kappa=\frac{1}{\rho}.
\end{equation*}
For $g\in C^{1,\delta}(\mathbb{S}^1)$, the Lax--Milgram theorem \cite[Chapter 5]{brezis} yields a unique periodic weak solution $f\in H^1(\mathbb{S}^1)$ of $f''-f=g$.
Testing with $f$, we obtain
\begin{equation*}
\int_0^{2\pi}(|f'|^2+f^2)\mathrm{d}\theta=-\int_0^{2\pi}gf\mathrm{d}\theta
\leq\|g\|_{L^2}\|f\|_{L^2}.
\end{equation*}
Since $f''=f+g$, the Sobolev embedding theorem yields
\begin{equation*}
\|f\|_{C^1(\mathbb{S}^1)}\leq C\|f\|_{H^2(\mathbb{S}^1)}\leq C\|g\|_{L^2(\mathbb{S}^1)}.
\end{equation*}
Using $f''=f+g$ and $f'''=f'+g'$, we obtain $\|f\|_{C^{3,\delta}}\leq C\|g\|_{C^{1,\delta}}$. Thus $f\mapsto f''-f$ is an isomorphism from $C^{3,\delta}(\mathbb{S}^1)$ onto $C^{1,\delta}(\mathbb{S}^1)$. Since $q>0$ and $q\in C^{2,\delta}(\mathbb{S}^1)$, the operator
\begin{equation*}
f\mapsto\frac{q^2}{2}(f''-f):C^{3,\delta}(\mathbb{S}^1)\to C^{1,\delta}(\mathbb{S}^1)
\end{equation*}
is also an isomorphism.
Substituting \eqref{eq:dq} into \eqref{eq:dk} and using $\rho\kappa=1$, we obtain
\begin{equation*}
D\mathcal{K}(h)[f]=\frac{q^2}{2}(f''-f)+qq'f'-q\rho(\partial_\nu w)\circ x_h.
\end{equation*}
Since $L_v$ is invertible on $E$ by the proof of Lemma \ref{lem:shape}, the Schauder estimate \cite[Chapter 6]{gt} for \eqref{eq:w} gives
\begin{equation}\label{ee:est-w1}
\|w\|_{C^{2,\delta}(\overline{\Omega})}\leq C\|(qf)\circ x_h^{-1}\|_{C^{2,\delta}(\partial\Omega)}\leq C\|f\|_{C^{2,\delta}(\mathbb{S}^1)}.
\end{equation}
Since $q\in C^{2,\delta}(\mathbb{S}^1)$ and $\rho\in C^{1,\delta}(\mathbb{S}^1)$, \eqref{ee:est-w1} implies
\begin{equation}\label{ee:est-w2}
\|qq'f'-q\rho(\partial_\nu w)\circ x_h\|_{C^{1,\delta}(\mathbb{S}^1)}\leq C\|f\|_{C^{2,\delta}(\mathbb{S}^1)}.
\end{equation}
By \eqref{ee:est-w2} and the compact embedding of $C^{3,\delta}(\mathbb{S}^1)$ into $C^{2,\delta}(\mathbb{S}^1)$, the map
\begin{equation*}
f\mapsto qq'f'-q\rho(\partial_\nu w)\circ x_h
\end{equation*}
is compact from $C^{3,\delta}(\mathbb{S}^1)$ to $C^{1,\delta}(\mathbb{S}^1)$. Hence $D\mathcal{K}(h)$, regarded as an operator from $C^{3,\delta}(\mathbb{S}^1)$ to $C^{1,\delta}(\mathbb{S}^1)$, is Fredholm of index zero.

\vskip 0.1cm
Let $f\in C^{3,\delta}(\mathbb{S}^1)$ satisfy $D\mathcal{K}(h)[f]=0$. For small $|t|$, set
\begin{equation*}
\Omega_t=\Omega_{h+tf},\quad v_t=v_{h+tf},\quad k_t=\mathcal{K}(h+tf).
\end{equation*}
Let $\Phi_t=\Phi_{tf}=\mathrm{Id}+tV_f$ be the diffeomorphism of Lemma \ref{lem:shape}. Writing $k=\mathcal{K}(h)$, we have
\begin{equation}\label{eq:dkzero}
\left.\frac{\mathrm{d}}{\mathrm{d}t}k_t\right|_{t=0}=D\mathcal{K}(h)[f]=0\quad\text{in }C^{1,\delta}(\mathbb{S}^1).
\end{equation}
By \eqref{eq:conv}, $D^2(v^{-\beta})>0$ in $\Omega$. Let $M$ be the matrix in \eqref{eq:match} constructed from the fixed reference solution $v$ in $\Omega$. Set
\begin{equation*}
T_t=D^2v_t-M(v_t,\nabla v_t)\quad\text{on }\overline\Omega_t.
\end{equation*}
Then $T_0=0$ by \eqref{eq:m}.
By Lemma \ref{lem:field}, $\operatorname{tr}T_t=0$ in $\Omega_t$ and
\begin{equation*}
(\operatorname{div}T_t)_i=-\sum_{j,\ell=1}^2\partial_{b_\ell}M_{ij}(v_t,\nabla v_t)(T_t)_{\ell j},\quad i=1,2.
\end{equation*}
Set $z_t=v_t\circ\Phi_t$ and $A_t=(I+tDV_f)^{-1}$. By the proof of Lemma \ref{lem:shape}, $t\mapsto z_t$ is $C^1$ in $C^{2,\delta}(\overline{\Omega})$. The map $t\mapsto A_t$ is $C^1$ in $C^{1,\delta}(\overline{\Omega})$. By the chain rule,
\begin{equation*}
(\nabla v_t)\circ\Phi_t=A_t^{\mathsf{T}}\nabla z_t,\quad
D^2z_t=(D\Phi_t)^{\mathsf{T}}\big((D^2v_t)\circ\Phi_t\big)D\Phi_t
+\sum_{i=1}^2\big((\partial_i v_t)\circ\Phi_t\big)D^2(\Phi_t)_i.
\end{equation*}
Thus $(\nabla v_t)\circ\Phi_t$ and $(D^2v_t)\circ\Phi_t$ are differentiable at $t=0$ in $C(\overline\Omega)$.
By Lemma \ref{lem:field}, $M$ is $C^1$ near the values of $(v,\nabla v)$ on $\overline\Omega$. Differentiating $T_t\circ\Phi_t$ at $t=0$ and using $T_0=0$, we obtain $T\in C(\overline\Omega;\mathbb{R}^{2\times2})$ such that
\begin{equation}\label{eq:dt}
\Big\|\frac{T_t\circ\Phi_t}{t}-T \Big\|_{C(\overline\Omega)}\to0\quad\text{as }t\to0.
\end{equation}
By \eqref{eq:dt}, $T$ is symmetric and $\operatorname{tr}T=0$.
Fix a compact set $K\subset\Omega$. For small $|t|$, we have $K\subset\Omega_t$ and $\Phi_t^{-1}\to\mathrm{Id}$ uniformly on $K$. By \eqref{eq:dt}, the convergence $z_t\to v$ in $C^2(\overline{\Omega})$ and the chain rule,
\begin{equation*}
\frac{T_t}{t}\to T,\quad (v_t,\nabla v_t)\to(v,\nabla v)\quad\text{uniformly on }K.
\end{equation*}
Fix $i\in\{1,2\}$ and $\xi\in C_c^\infty(\Omega)$. Testing the equation for $(\operatorname{div}T_t)_i$ with $\xi/t$ and passing to the limit on $\operatorname{supp}\xi$, we obtain
\begin{equation*}
\sum_{j=1}^2\int_\Omega T_{ij}\partial_{x_j}\xi\mathrm{d}x
=\sum_{j,\ell=1}^2\int_\Omega\partial_{b_\ell}M_{ij}(v,\nabla v)T_{\ell j}\xi\mathrm{d}x.
\end{equation*}
Hence
\begin{equation*}
(\operatorname{div}T)_i=-\sum_{j,\ell=1}^2\partial_{b_\ell}M_{ij}(v,\nabla v)T_{\ell j}\quad\text{in }\mathcal{D}'(\Omega).
\end{equation*}
Since $D_bM(v,\nabla v)$ is bounded on $\overline{\Omega}$ by Lemma \ref{lem:field}, we have $\operatorname{div}T\in L^\infty(\Omega;\mathbb{R}^2)$ and $|\operatorname{div}T|\leq C|T|$ a.e. in $\Omega$.

\vskip 0.1cm
At $x_{h+tf}(\theta)=\Phi_t(x_h(\theta))$, the normal and tangent are $n(\theta)$ and $\tau(\theta)$, independent of $t$.
Recall that $\gamma=\frac{1}{\alpha+1}$ and $m=2-\gamma$ where $\alpha=\frac{1}{\beta}=\frac{2}{p-1}$ for $p>1$.
At $x_{h+tf}(\theta)$, \eqref{eq:bd} gives
\begin{equation}\label{eq:bd-var}
T_t(n,n)=\frac{q_{h+tf}^{1+\gamma}}{2}\left(\frac{q_{h+tf}^m}{k_t}-\frac{q^m}{k}\right),\quad
T_t(n,\tau)=-\frac{q_{h+tf}^{1+\gamma}}{2m}\left(\frac{(q_{h+tf}^m)'}{k_t}-\frac{(q^m)'}{k}\right).
\end{equation}
Define $H$ on $\partial\Omega$ by $H(x_h(\theta))=mq(\theta)^{m-1}(Dq_h[f])(\theta)$. By \eqref{eq:dkzero},
\begin{equation*}
\left.\frac{\mathrm{d}}{\mathrm{d}t}\frac{q_{h+tf}^m}{k_t}\right|_{t=0}=\frac{H\circ x_h}{k},\quad
\left.\frac{\mathrm{d}}{\mathrm{d}t}\frac{(q_{h+tf}^m)'}{k_t}\right|_{t=0}=\frac{(H\circ x_h)'}{k}.
\end{equation*}
By \eqref{eq:k} and $x_h'=\tau/\kappa$ in \eqref{eq:param},
\begin{equation*}
2\kappa k=q^2,\quad (H\circ x_h)'(\theta)=\frac{\partial_\tau H(x_h(\theta))}{\kappa(\theta)}.
\end{equation*}
Differentiating \eqref{eq:bd-var} at $t=0$ and using \eqref{eq:dt}, we obtain at $x_h(\theta)$
\begin{equation*}
\begin{aligned}
T(\nu,\nu)&=\frac{q^{1+\gamma}}{2k}H=\kappa q^{\gamma-1}H,\\
T(\nu,\tau)&=-\frac{q^{1+\gamma}}{2mk}(H\circ x_h)'(\theta)=-\frac{q^{\gamma-1}}{m}\partial_\tau H.
\end{aligned}
\end{equation*}
Lemma \ref{lem:tensor} gives $T=0$ in $\Omega$. By continuity of $T$ on $\overline{\Omega}$, $T=0$ on $\partial\Omega$ as well. Since $\kappa,q>0$, the identity for $T(\nu,\nu)$ yields $H=0$ on $\partial\Omega$. Hence $Dq_h[f]=0$ on $\mathbb{S}^1$, and \eqref{eq:dk} gives
\begin{equation*}
0=D\mathcal{K}(h)[f]=\frac{q^2}{2}(f+f'').
\end{equation*}
Consequently $f=b\cdot n$ for some $b\in\mathbb{R}^2$. Conversely, $h+tb\cdot n$ is the support function of $\Omega_h+tb$, so $\mathcal{K}(h+tb\cdot n)=\mathcal{K}(h)$. Differentiating at $t=0$, we obtain $D\mathcal{K}(h)[b\cdot n]=0$. Hence
\begin{equation*}
\ker\left(D\mathcal{K}(h):C^{3,\delta}(\mathbb{S}^1)\to C^{1,\delta}(\mathbb{S}^1)\right)=\operatorname{span}\{\cos\theta,\sin\theta\}.
\end{equation*}
Since $D\mathcal{K}(h)$ has index zero on $C^{3,\delta}(\mathbb{S})$, its range has codimension two in $C^{1,\delta}(\mathbb{S})$. For every $f\in C^{3,\delta}(\mathbb{S}^1)$, we have $\mathcal{K}(h+tf)\in Y$ for small $|t|$. Hence
\begin{equation*}
0=\left.\frac{\mathrm{d}}{\mathrm{d}t}\int_0^{2\pi}\mathcal{K}(h+tf)n\mathrm{d}\theta\right|_{t=0}
=\int_0^{2\pi}D\mathcal{K}(h)[f]n\mathrm{d}\theta.
\end{equation*}
Thus the range is contained in $Y$. Since $\int_0^{2\pi}n\otimes n\mathrm{d}\theta=\pi I$, the space $Y$ also has codimension two. Hence the range of $D\mathcal{K}(h)$ is $Y$.

\vskip 0.1cm
Every $f\in C^{3,\delta}(\mathbb{S})$ has the unique decomposition
\begin{equation*}
f=f_0+b\cdot n,\quad f_0\in X,\quad b=\frac{1}{\pi}\int_0^{2\pi}fn\mathrm{d}\theta.
\end{equation*}
Since the kernel consists precisely of the functions $b\cdot n$, the restriction of $D\mathcal{K}(h)$ to $X$ is a bijection from $X$ to $Y$. Its inverse is bounded by the bounded inverse theorem \cite[Chapter 2]{brezis}.
\end{proof}

For $h\in\mathcal{A}$, define the ground-state energy by
\begin{equation}\label{eq:nehari}
\mathcal E(h):=\frac{p-1}{2(p+1)}S_p(\Omega_h)^{\frac{p+1}{p-1}} =\frac{p-1}{2(p+1)} \int_{\Omega_h}v_h^{p+1}\mathrm{d}x,
\end{equation}
where the second equality follows from \eqref{eq:least}.

\begin{lemma}\label{lem:energy}
The functional $\mathcal E : \mathcal A\to\mathbb{R}$ is $C^1$. For $h\in\mathcal A$ and $f\in X$,
\begin{equation}\label{eq:dc}
D\mathcal{E}(h)[f]=-\frac{1}{2}\int_{\partial\Omega_h}(\partial_\nu v_h)^2(f\circ x_h^{-1})\mathrm{d}s=-\int_0^{2\pi}\mathcal{K}(h)f\mathrm{d}\theta.
\end{equation}
Moreover,
\begin{equation}\label{eq:sc}
2\alpha \mathcal E(h) =\int_0^{2\pi}h\mathcal{K}(h)\mathrm{d}\theta.
\end{equation}
\end{lemma}

\begin{proof}
Fix $h\in\mathcal{A}$ and write $\Omega=\Omega_h$ and $v=v_h$. For $g\in X$ sufficiently small, set $z_g=v_{h+g}\circ\Phi_g$. By \eqref{eq:nehari} and the change of variables $y=\Phi_g(x)$,
\begin{equation*}
\mathcal{E}(h+g)=\frac{p-1}{2(p+1)}\int_\Omega z_g^{p+1}\det D\Phi_g\mathrm{d}x.
\end{equation*}
By the proof of Lemma \ref{lem:shape}, $g\mapsto z_g$ is smooth with values in $C^{2,\delta}_0(\overline{\Omega})$. Since $g\mapsto\det D\Phi_g$ is smooth and $s\mapsto(s_+)^{p+1}$ is $C^1$, we have $\mathcal{E}\in C^1(\mathcal{A})$.

For $f\in X$, write $\zeta_f=\left.\frac{\mathrm{d}}{\mathrm{d}t}z_{tf}\right|_{t=0}$ and $w=\zeta_f-\nabla v\cdot V_f$. Differentiating the formula for $\mathcal{E}(h+tf)$ at $t=0$ and using
\begin{equation*}
\left.\frac{\mathrm{d}}{\mathrm{d}t}\det(I+tDV_f)\right|_{t=0}=\operatorname{div}V_f,
\end{equation*}
we obtain
\begin{equation}\label{eq:3.4-1}
\begin{aligned}
\frac{2}{p-1}D\mathcal E(h)[f]
&=\int_\Omega v^p\zeta_f\mathrm{d}x+\frac{1}{p+1}\int_\Omega v^{p+1}\operatorname{div}V_f\mathrm{d}x\\
&=\int_\Omega v^p w\mathrm{d}x+\frac{1}{p+1}\int_\Omega\operatorname{div}(v^{p+1}V_f)\mathrm{d}x=\int_\Omega v^p w\mathrm{d}x.
\end{aligned}
\end{equation}
By \eqref{eq:w}, $L_vw=0$ in $\Omega$ and $w(x_h(\theta))=q_h(\theta)f(\theta)$. Since $L_vv=-(p-1)v^p$ in $\Omega$ and $v=0$ on $\partial\Omega$, Green's identity yields
\begin{equation}\label{eq:3.4-2}
\begin{aligned}
(p-1)\int_\Omega v^p w\mathrm{d}x
&=\int_\Omega(vL_vw-wL_vv)\mathrm{d}x=\int_{\partial\Omega}w\partial_\nu v\mathrm{d}s =-\int_{\partial\Omega}(\partial_\nu v)^2(f\circ x_h^{-1})\mathrm{d}s.
\end{aligned}
\end{equation}
Noting that $\mathrm{d}s=(h+h'')\mathrm{d}\theta$, we obtain \eqref{eq:dc} from \eqref{eq:3.4-1}-\eqref{eq:3.4-2} (see also \cite[Lemma 1]{cfr}).

For $s>0$, the domain with support function $sh$ is $s\Omega_h$.
By the change of variables in \eqref{eq:sob} and the definition of $\mathcal E$ in \eqref{eq:nehari}, we have
\begin{equation}\label{eq:sscale}
S_p(s\Omega_h)=s^{-\frac{4}{p+1}}S_p(\Omega_h),\quad
\mathcal E(sh)=s^{-2\alpha}\mathcal E(h).
\end{equation}
Differentiating $\mathcal E(sh)$ at $s=1$ and taking $f=h$ in \eqref{eq:dc}, we obtain
\begin{equation*}
-2\alpha \mathcal E(h)=D\mathcal E(h)[h]=-\int_0^{2\pi}h(\theta)\mathcal{K}(h)(\theta)\mathrm{d}\theta,
\end{equation*}
which gives \eqref{eq:sc}.
\end{proof}

Let $v_1$ be the positive ground state of problem  \eqref{eq:le} in $B_1(0)$. By rotation invariance and Proposition \ref{prop:bf}, $v_1$ is radial. Thus $\mathcal{K}(1)$ is a positive constant. Set $v_R(x)=R^{-\alpha}v_1(x/R)$. Since $\alpha p=\alpha+2$, $v_R$ solves \eqref{eq:le} in $B_R(0)$. By \eqref{eq:least} and \eqref{eq:sscale},
\begin{equation*}
\int_{B_R(0)}v_R^{p+1}\mathrm{d}x
=R^{-2\alpha}\int_{B_1(0)}v_1^{p+1}\mathrm{d}x
=S_p(B_R(0))^{\frac{p+1}{p-1}}.
\end{equation*}
Thus $v_R$ is the positive ground state of \eqref{eq:le} in $B_R(0)$. Since
\begin{equation*}
q_{v_R}(Rn(\theta))=R^{-\alpha-1}q_{v_1}(n(\theta)),\quad \kappa=\frac{1}{R}\quad\text{on }\partial B_R(0),
\end{equation*}
\eqref{eq:map} gives $\mathcal{K}(R)=R^{-2\alpha-1}\mathcal{K}(1)$. Set $k_0=\max_{\mathbb{S}}k$ and choose $R_0>0$ such that $\mathcal{K}(R_0)=k_0$. For $0\leq t\leq1$, set
\begin{equation}\label{eq:def-kt1}
k_t=(1-t)k_0+tk.
\end{equation}
Then $0<\min_{\mathbb{S}}k\leq k_t\leq k_0$. By \eqref{eq:target} and $\int_0^{2\pi}n(\theta)\mathrm{d}\theta=0$, we have $k_t\in Y$.

\begin{lemma}\label{lem:bound}
Suppose that $t\mapsto h_t\in\mathcal{A}$ is a continuously differentiable curve on $[0,T)$ for some $T\in(0,1]$ satisfying $h_0=R_0$ and
\begin{equation}\label{eq:path}
\mathcal{K}(h_t)=k_t\quad\text{for all }t\in[0,T).
\end{equation}
Then, for each $\delta_1\in(\delta,1)$, we have
\begin{equation}\label{eq:hbound}
\|h_t\|_{C^{3,\delta_1}(\mathbb{S}^1)}\leq C,\quad h_t+h_t''\geq c\quad\text{on }\mathbb{S}^1, ~~\mbox{for all}~~~t\in[0,T),
\end{equation}
where $c,C$ are positive constants depending only on $p$, $\delta_1$ and $k$.
\end{lemma}

\begin{proof}
For $0\leq t<T$, write
\begin{equation*}
\Omega_t=\Omega_{h_t},\quad v_t=v_{h_t},\quad q_t=q_{h_t},\quad \rho_t=h_t+h_t''.
\end{equation*}
Since $h_t\in\mathcal{A}$, the normalization and \eqref{eq:support} give $h_t>0$ on $\mathbb{S}$.
Using $k_t\geq\min_{\mathbb{S}^1}k$ and $h_t>0$, \eqref{eq:sc} and \eqref{eq:path} give
\begin{equation*}
2\alpha \mathcal E(h_t)=\int_0^{2\pi}h_t k_t\mathrm{d}\theta\geq \Big(\min_{\mathbb{S}^1}k \Big) \int_0^{2\pi}h_t\mathrm{d}\theta.
\end{equation*}
By \eqref{eq:dc} and \eqref{eq:path},
\begin{equation*}
\int_0^{2\pi}k_t \partial_t h_t\mathrm{d}\theta=-\frac{\mathrm{d}}{\mathrm{d}t}\mathcal{E}(h_t).
\end{equation*}
Differentiating \eqref{eq:sc} and using $\partial_t k_t=k-k_0$, $h_t>0$ and $k\leq k_0$, we obtain
\begin{equation}\label{eq:de<0}
(2\alpha+1)\frac{\mathrm{d}}{\mathrm{d}t}\mathcal{E}(h_t)=\int_0^{2\pi}h_t(k-k_0)\mathrm{d}\theta\leq0,
\end{equation}
which gives
\begin{equation}\label{eq:cbd}
\mathcal E(h_t)\leq \mathcal E(h_0) \quad \text{and} \quad \int_0^{2\pi}h_t\mathrm{d}\theta\leq\frac{2\alpha \mathcal E(h_0)}{\min_{\mathbb{S}^1}k}.
\end{equation}
For $x\in\Omega_t$, we have $(x\cdot n(\theta))_+\leq h_t(\theta)$. Hence
\begin{equation*}
2|x|=\int_0^{2\pi}(x\cdot n(\theta))_+\mathrm{d}\theta
\leq\int_0^{2\pi}h_t(\theta)\mathrm{d}\theta\leq C.
\end{equation*}
Thus $\Omega_t\subset B_R(0)$ for some $R>0$ independent of $t$. Write $B_R=B_R(0)$.
By \eqref{eq:energy}, \eqref{eq:nehari} and \eqref{eq:cbd}, we have
\begin{equation*}
\int_{\Omega_t}|\nabla v_t|^2\mathrm{d}x=\int_{\Omega_t} v^{p+1}_t\mathrm{d}x=\frac{2(p+1)}{p-1} \mathcal E(h_t)\leq C.
\end{equation*}
The zero extension of $v_t$ belongs to $H^1_0(B_R)$. By the Sobolev and Poincar\'e inequalities on $B_R$,
\begin{equation*}
\|v^p_t\|_{L^2(\Omega_t)}=\|v_t\|_{L^{2p}(\Omega_t)}^p\leq C\|\nabla v_t\|_{L^2(\Omega_t)}^p\leq C.
\end{equation*}
Let $z_t\in W^{2,2}(B_R)\cap H^1_0(B_R)$ solve
\begin{equation*}
\begin{cases}
-\Delta z_t=v_t^p\mathbf 1_{\Omega_t} & \text{in }B_R,\\
z_t=0 & \text{on }\partial B_R.
\end{cases}
\end{equation*}
The maximum principle gives $z_t\geq0$ in $B_R$. Since $z_t-v_t$ is harmonic in $\Omega_t$ and $z_t-v_t=z_t\geq0$ on $\partial\Omega_t$, we have $0\leq v_t\leq z_t$ in $\Omega_t$. The $W^{2,2}$ estimate \cite[Chapter 9]{gt} on $B_R$ and the Sobolev embedding theorem yield
\begin{equation*}
\|v_t\|_{L^\infty(\Omega_t)}\leq\|z_t\|_{L^\infty(B_R)}\leq C\|v^p_t\|_{L^2(\Omega_t)}\leq C.
\end{equation*}
On the other hand, by \eqref{eq:energy} and the Poincar\'e inequality on $B_R$, we have
\begin{equation*}
0<\int_{\Omega_t}|\nabla v_t|^2\mathrm{d}x
=\int_{\Omega_t} v_t^{p+1}\mathrm{d}x \leq C_R\|v_t\|_{L^\infty(\Omega_t)}^{p-1}\int_{\Omega_t}|\nabla v_t|^2\mathrm{d}x.
\end{equation*}
Hence $\displaystyle\max_{\overline{\Omega_t}}v_t\geq C_R^{-\frac{1}{p-1}}>0$.

\vskip 0.1cm
We next prove a uniform lower bound for $\rho_t$. Choose $C_0>0$ such that $\|v_t^p\|_{L^\infty(\Omega_t)}\leq C_0$ for all $0\leq t<T$. Fix $y\in\partial\Omega_t$ and set $d_y(x)=(y-x)\cdot\nu(y)$. By convexity and $\Omega_t\subset B_R$, we have $0\leq d_y\leq 2R$ on $\overline{\Omega_t}$. Set
\begin{equation*}
b_y(x)=\frac{C_0}{2}d_y(x)(2R-d_y(x)).
\end{equation*}
Since $\nabla d_y=-\nu(y)$ and $\Delta d_y=0$, we obtain
\begin{equation*}
\begin{cases}
-\Delta(b_y-v_t)=C_0-v_t^p\geq0 &\text{in }\Omega_t,\\
b_y-v_t=b_y\geq0 &\text{on }\partial\Omega_t.
\end{cases}
\end{equation*}
The maximum principle yields
\begin{equation}\label{e:vtc0}
0\leq v_t(x)\leq b_y(x)\leq C_0R d_y(x)\quad\text{for }x\in\Omega_t.
\end{equation}

For each $x\in\Omega_t$, choose a nearest boundary point $y$. Then $x-y$ is normal to $\partial\Omega_t$ at $y$, and convexity gives $(x-y)\cdot\nu(y)\leq0$. Hence $d_y(x)=|x-y|$, and \eqref{e:vtc0} gives
\begin{equation*}
v_t(x)\leq C_0R\operatorname{dist}(x,\partial\Omega_t).
\end{equation*}
At $y=x_{h_t}(\theta)$, we have $d_y(y-s\nu(y))=s$. Applying \eqref{e:vtc0} at $y-s\nu(y)$ and letting $s\downarrow0$, we obtain
\begin{equation*}
q_t(\theta)=-\partial_\nu v_t(y)=\lim_{s\downarrow0}\frac{v_t(y-s\nu(y))}{s}\leq C_0R.
\end{equation*}
Since $k_t=(1-t)k_0+tk$ and $k_0=\max_{\mathbb{S}^1}k$, we have $k_t\geq\min_{\mathbb{S}^1}k$. By \eqref{eq:map} and \eqref{eq:path},
\begin{equation*}
\rho_t=\frac{2k_t}{q_t^2}\geq r_0,\quad r_0=\frac{2\min_{\mathbb{S}^1}k}{C_0^2R^2}>0.
\end{equation*}
It follows from \eqref{eq:disk} that
\begin{equation*}
B_{r_0}(y-r_0\nu(y))\subset\Omega_t\quad\text{for every }y\in\partial\Omega_t.
\end{equation*}

To bound $q_t$ from below, take $x_*\in{\Omega_t}$ satisfying $v_t(x_*)=\max_{\overline{\Omega_t}}v_t\geq C_R^{-\frac{1}{p-1}}>0$.
Since $v_t(x_*)\leq C_0R\operatorname{dist}(x_*,\partial{\Omega_t})$, we obtain
\begin{equation*}
v_t(x_*)\geq c>0,\quad \operatorname{dist}(x_*,\partial{\Omega_t})\geq c_0>0.
\end{equation*}
Fix $y\in\partial\Omega_t$ and set
\begin{equation*}
a=y-r_0\nu(y),\quad r_1=\frac{1}{4}\min\left\{c_0,\frac{r_0}{2}\right\}.
\end{equation*}
For every $z\in\overline B_{r_0/2}(a)$, both $B_{4r_1}(x_*)$ and $B_{4r_1}(z)$ lie inside ${\Omega_t}$.
By the convexity of $\Omega_t$, we get
\begin{equation*}
B_{4r_1}(x)\subset{\Omega_t}\quad\text{for every }x\in[x_*,z].
\end{equation*}
Choose points $x_0=x_*,x_1,\cdots,x_\ell=z$ on $[x_*,z]$ such that
\begin{equation*}
|x_i-x_{i-1}|\leq\frac{r_1}{2}\quad\text{for }1\leq i\leq\ell,\quad \ell\leq1+\left\lceil\frac{4R}{r_1}\right\rceil.
\end{equation*}
Since $\|v_t^{p-1}\|_{L^\infty(\Omega_t)}\leq C$, the Harnack inequality \cite[Chapter 8]{gt} for $-\Delta v_t=v_t^{p-1}v_t$ in $B_{2r_1}(x_{i-1})$ yields
\begin{equation*}
v_t(x_{i-1})\leq\sup_{B_{r_1}(x_{i-1})}v_t
\leq C\inf_{B_{r_1}(x_{i-1})}v_t \leq C v_t(x_i),
\end{equation*}
where the constant $C>0$ is independent of $t$, $y$ and $z$.
Thus,
\begin{equation*}
v_t(z)\geq C^{-\ell}v_t(x_*)\geq c_1>0\quad\text{for every }z\in\overline B_{r_0/2}(a).
\end{equation*}
On the annulus $r_0/2<|x-a|<r_0$, set
\begin{equation*}
w(x)=\frac{c_1}{\log 2}\log\frac{r_0}{|x-a|}.
\end{equation*}
Since $w=c_1$ on $\partial B_{r_0/2}(a)$ and $w=0$ on $\partial B_{r_0}(a)$, we have
\begin{equation*}
\begin{cases}
-\Delta(v_t-w)=v_t^p\geq0&\text{in }B_{r_0}(a)\setminus\overline B_{\frac{r_0}{2}}(a),\\[1mm]
v_t-w\geq0&\text{on }\partial B_{r_0}(a)\cup\partial B_{\frac{r_0}{2}}(a).
\end{cases}
\end{equation*}
Write $y=x_{h_t}(\theta)$. For $0<s<r_0/2$, the maximum principle gives
\begin{equation*}
v_t(y-s\nu(y))\geq\frac{c_1}{\log 2}\log\frac{r_0}{r_0-s}.
\end{equation*}
Dividing by $s$ and letting $s\to 0^{+}$, we obtain
\begin{equation*}
q_t(\theta)=-\partial_\nu v_t(y)\geq\frac{c_1}{\log 2}\lim_{s\to 0^{+}}\frac{1}{s}\log\frac{r_0}{r_0-s}=\frac{c_1}{r_0\log 2}.
\end{equation*}
Since $q_t\leq C_0R$ and $\rho_t=2k_t/q_t^2$,
\begin{equation}\label{eq:rho}
0<c\leq q_t\leq C,\quad 0<c\leq\rho_t=\frac{2k_t}{q_t^2}\leq C.
\end{equation}

Fix $y=x_{h_t}(\theta_0)\in\partial\Omega_t$. Use coordinates centered at $y$ with axes along $\tau(\theta_0)$ and $-n(\theta_0)$. Write $(x_1(\vartheta),x_2(\vartheta))$ for the coordinates of $x_{h_t}(\theta_0+\vartheta)$. By \eqref{eq:param},
\begin{equation*}
x_1'(\vartheta)=\rho_t(\theta_0+\vartheta)\cos\vartheta,\quad
x_2'(\vartheta)=\rho_t(\theta_0+\vartheta)\sin\vartheta.
\end{equation*}
For $|\vartheta|\leq\pi/4$, the bound $\rho_t\geq r_0$ gives
\begin{equation*}
x_1'(\vartheta)\geq\frac{r_0}{\sqrt{2}},\quad
\pm x_1(\pm\pi/4)=\int_0^{\pi/4}\rho_t(\theta_0\pm s)\cos s\mathrm{d}s\geq\frac{r_0}{\sqrt{2}}.
\end{equation*}
Therefore, the boundary arc of $\Omega_t$ can be described by a graph $x_2=g(x_1)$ for $|x_1|<\frac{r_0}{2}$. Since $g'=x_2'/x_1'=\tan\vartheta$,
\begin{equation}\label{eq:graph}
\frac{\mathrm{d}\vartheta}{\mathrm{d}x_1}=\frac{1}{\rho_t(\theta_0+\vartheta)\cos\vartheta},\quad
g''=\frac{1}{\rho_t(\theta_0+\vartheta)\cos^3\vartheta}.
\end{equation}
Fix $|x_1|<\frac{r_0}{2}$.
Since $g(0)=0$ and $|g'|\leq1$, we have $g(x_1)\leq|x_1|<r_0$.
By \eqref{eq:disk}, $(x_1,r_0)\in B_{r_0}((0,r_0))\subset{\Omega_t}$.
By convexity of $\Omega_t$, the vertical section of $\Omega_t$ at $x_1$ is an interval containing $(x_1,r_0)$. Since $g(x_1)<r_0$ and $(x_1,g(x_1))\in\partial\Omega_t$, its lower endpoint is $(x_1,g(x_1))$. Hence
\begin{equation*}
{\Omega_t}\cap B_{{r_0}/{4}}(0)=\big\{(x_1,x_2)\in B_{{r_0}/{4}}(0):x_2>g(x_1)\big\}.
\end{equation*}
We also get from \eqref{eq:rho} and \eqref{eq:graph} that
\begin{equation*}
\|g\|_{C^{1,1}([-\frac{r_0}{4},\frac{r_0}{4}])}\leq C.
\end{equation*}
Fix $r>\frac{2}{1-\delta_1}$, so that $1-\frac{2}{r}>\delta_1$.
The boundary graphs have a common radius $r_0/4$ and uniformly bounded $C^{1,1}$ norms. The $W^{2,r}$ estimate \cite[Chapter 9]{gt} and the Sobolev embedding theorem yield
\begin{equation*}
\|v_t\|_{C^{1,\delta_1}(\overline{\Omega_t})}\leq C\left(\|v_t\|_{L^\infty({\Omega_t})}+\|v_t^p\|_{L^\infty({\Omega_t})}\right)\leq C.
\end{equation*}
Since $x_{h_t}'=\rho_t\tau$ and $\rho_t$ is uniformly bounded by \eqref{eq:rho}, the maps $x_{h_t}$ are uniformly Lipschitz. From $q_t(\theta)=-\nabla v_t(x_{h_t}(\theta))\cdot n(\theta)$, we obtain
\begin{equation*}
\|q_t\|_{C^{0,\delta_1}(\mathbb{S})}\leq C\|v_t\|_{C^{1,\delta_1}(\overline{\Omega_t})}\leq C.
\end{equation*}
Since $k_t=(1-t)k_0+tk$ and $k\in C^2(\mathbb{S}^1)$, the functions $k_t$ are uniformly bounded in $C^{1,\delta_1}(\mathbb{S}^1)$. By $\rho_t=2k_t/q_t^2$, the $C^{0,\delta_1}$ bound for $q_t$, and \eqref{eq:rho},
\begin{equation}\label{ee:3-1}
\|\rho_t\|_{C^{0,\delta_1}(\mathbb{S})}+\|\rho_t^{-1}\|_{C^{0,\delta_1}(\mathbb{S})}\leq C.
\end{equation}
By \eqref{eq:graph}, $\vartheta(x_1)$ is uniformly Lipschitz and $\cos\vartheta\geq1/\sqrt{2}$. Substituting \eqref{ee:3-1} into the formula for $g''$, we obtain
\begin{equation}
\|g\|_{C^{2,\delta_1}([-r_0/4,r_0/4])}\leq C. \label{ee:3-3}
\end{equation}
Choose $M>0$ such that $\|v_t\|_{L^\infty(\Omega_t)}\leq M$ for all $0\leq t<T$. Since $s\mapsto s^p$ is Lipschitz on $[0,M]$, the uniform $C^{1,\delta_1}$ bound for $v_t$ gives $\|v_t^p\|_{C^{0,\delta_1}(\overline{\Omega_t})}\leq C$. By \eqref{ee:3-3} and the Schauder estimate \cite[Chapter 6]{gt},
\begin{equation}\label{ee:3-4}
\|v_t\|_{C^{2,\delta_1}(\overline{\Omega_t})}\leq C\left(\|v_t\|_{L^\infty({\Omega_t})}+\|v_t^p\|_{C^{0,\delta_1}(\overline{\Omega_t})}\right)\leq C.
\end{equation}
Differentiating $q_t(\theta)=-\nabla v_t(x_{h_t}(\theta))\cdot n(\theta)$ and using $\nabla v_t\cdot\tau=0$ on $\partial{\Omega_t}$, we obtain
\begin{equation*}
q_t'(\theta)=-\rho_t(\theta)D^2v_t(x_{h_t}(\theta))[n(\theta),\tau(\theta)].
\end{equation*}
By \eqref{ee:3-1}, \eqref{ee:3-4} and the uniform Lipschitz bound for $x_{h_t}$, the functions $q_t'$ are uniformly bounded in $C^{0,\delta_1}(\mathbb{S}^1)$. Using $\rho_t=2k_t/q_t^2$ and $q_t\geq c$, we obtain
\begin{equation*}
\|q_t\|_{C^{1,\delta_1}(\mathbb{S}^1)}\leq C,\quad
\|\rho_t\|_{C^{1,\delta_1}(\mathbb{S}^1)}= \Big\|\frac{2k_t}{q_t^2} \Big\|_{C^{1,\delta_1}(\mathbb{S}^1)}\leq C.
\end{equation*}
Since $\Omega_t\subset B_R(0)$, the identities $h_t=x_{h_t}\cdot n$ and $h_t'=x_{h_t}\cdot\tau$ give
\begin{equation*}
\|h_t\|_{C^1(\mathbb{S}^1)}\leq 2R.
\end{equation*}
By \eqref{eq:rho}, $h_t''=\rho_t-h_t$ is uniformly bounded, so $h_t'$ is uniformly Lipschitz. Using $h_t'''=\rho_t'-h_t'$ and the uniform $C^{1,\delta_1}$ bound for $\rho_t$, we obtain
\begin{equation*}
\|h_t\|_{C^{3,\delta_1}(\mathbb{S}^1)}\leq C.
\end{equation*}
The lower bound $h_t+h_t''=\rho_t\geq c$ follows from \eqref{eq:rho}. This proves \eqref{eq:hbound}.
\end{proof}

Let $h_1,h_2\in\mathcal{A}$ satisfy $\mathcal{K}(h_1)=\mathcal{K}(h_2)$. By \eqref{eq:map} and Lemma \ref{lem:conv},
\begin{equation*}
k_{v_{h_1}}=k_{v_{h_2}} \quad\text{on }\mathbb{S}^1,\quad
D^2(v_{h_2}^{-\beta})>0 \quad\text{in }\Omega_{h_2}.
\end{equation*}
Proposition \ref{prop:comp} then gives $\Omega_{h_1}=a+\Omega_{h_2}$ for some $a\in\mathbb{R}^2$. Writing $x=a+y$ in \eqref{def-h}, we obtain
\begin{equation*}
h_1(\theta)=\max_{y\in\overline{\Omega_{h_2}}}(a+y)\cdot n(\theta)=a\cdot n(\theta)+h_2(\theta).
\end{equation*}
Since $h_1,h_2\in X$, the normalization in \eqref{eq:def-x+y} gives
\begin{equation*}
0=\int_0^{2\pi}(h_1-h_2)n\mathrm{d}\theta=\int_0^{2\pi}(a\cdot n)n\mathrm{d}\theta=\pi a.
\end{equation*}
Thus $a=0$ and $h_1=h_2$, so $\mathcal{K}$ is injective on $\mathcal{A}$.

\begin{proposition}\label{prop:real}
Let $k\in C^2(\mathbb{S}^1)$ satisfy \eqref{eq:target}. Then there exists a unique $h\in\mathcal{A}$ such that
\begin{equation*}
\mathcal{K}(h)=k.
\end{equation*}
\end{proposition}

\begin{proof}
By Lemmas \ref{lem:shape} and \ref{lem:inv}, the inverse function theorem yields a smooth local inverse of $\mathcal{K}$ near $\mathcal{K}(h)$ for every $h\in\mathcal{A}$. Hence $\mathcal{K}(\mathcal{A})$ is open in $Y$. Since $\mathcal{K}$ is injective on $\mathcal{A}$, the local inverses define a smooth inverse $\mathcal{K}^{-1}$ on $\mathcal{K}(\mathcal{A})$.

\vskip 0.1cm
Let $k_t=(1-t)k_0+tk$ be as in \eqref{eq:def-kt1}, where $k_0=\max_{\mathbb{S}^1}k=\mathcal{K}(R_0)$. Set
\begin{equation*}
T=\sup\big\{t\in[0,1]\mid k_s\in\mathcal{K}(\mathcal{A})\ \text{for every }s\in[0,t]\big\}.
\end{equation*}
Since $k_0$ belongs to the open set $\mathcal{K}(\mathcal{A})$ and $t\mapsto k_t$ is continuous in $Y$, we have $T>0$. For $0\leq t<T$, the curve $h_t=\mathcal{K}^{-1}(k_t)$ is $C^1$ and satisfies \eqref{eq:path} with $h_0=R_0$.

\vskip 0.1cm
Fix $\delta<\delta_1<1$. By Lemma \ref{lem:bound}, there exist $c,C>0$ such that for every $0\leq t<T$,
\begin{equation*}
\|h_t\|_{C^{3,\delta_1}(\mathbb{S}^1)}\leq C,\quad h_t+h_t''\geq c\quad\text{on }\mathbb{S}^1.
\end{equation*}
Choose $t_j\uparrow T$ with $t_j<T$. By the compact embedding of $C^{3,\delta_1}(\mathbb{S}^1)$ into $C^{3,\delta}(\mathbb{S}^1)$, after passing to a subsequence, we have
\begin{equation*}
h_{t_j}\to h_*\quad\text{in }C^{3,\delta}(\mathbb{S}^1).
\end{equation*}
Since $h_{t_j}\in X$ and $h_{t_j}+h_{t_j}''\geq c$, passing to the limit gives
\begin{equation*}
\begin{aligned}
\int_0^{2\pi}h_*(\theta)n(\theta)\mathrm{d}\theta&=\lim_{j\to\infty}\int_0^{2\pi}h_{t_j}(\theta)n(\theta)\mathrm{d}\theta=0,
\end{aligned}
\end{equation*}
and
\begin{equation*}
\begin{aligned}
h_*+h_*''&=\lim_{j\to\infty}(h_{t_j}+h_{t_j}'')\geq c\quad\text{on }\mathbb{S}^1.
\end{aligned}
\end{equation*}
Thus $h_*\in\mathcal{A}$. By continuity of $\mathcal{K}$,
\begin{equation*}
\mathcal{K}(h_*)=\lim_{j\to\infty}\mathcal{K}(h_{t_j})=\lim_{j\to\infty}k_{t_j}=k_T.
\end{equation*}

If $T<1$, then $k_T=\mathcal{K}(h_*)$ belongs to the open set $\mathcal{K}(\mathcal{A})$. By continuity of $t\mapsto k_t$, there exists $0<\varepsilon<1-T$ such that
\begin{equation*}
k_t\in\mathcal{K}(\mathcal{A})\quad\text{for every }t\in[T,T+\varepsilon].
\end{equation*}
Since $k_t\in\mathcal{K}(\mathcal{A})$ also for $0\leq t<T$, this contradicts the definition of $T$. Hence $T=1$ and $\mathcal{K}(h_*)=k$. Uniqueness follows from the injectivity of $\mathcal{K}$ on $\mathcal{A}$.
\end{proof}

\section{Proof of the main theorem}\label{sec:proof}

In this section, we prove Theorem \ref{thm:main} and Corollary \ref{cor:conv}. We first consider the bounded convex domains with positive boundary curvature. If the curvature vanishes somewhere, we approximate both the domain and the given positive solution.

\begin{proposition}\label{prop:pos}
Let $u$ be a positive solution of problem  \eqref{eq:le} in a bounded smooth convex domain $\Omega\subset\mathbb{R}^2$ with positive boundary curvature. Then $u$ is a ground state.
\end{proposition}

\begin{proof}
By \eqref{eq:reg}, $q_u=-\partial_\nu u$ is positive and belongs to $C^{2,\delta}(\partial\Omega)$. Since $\partial\Omega$ is smooth and $\kappa>0$, the boundary point $x(\theta)$ depends smoothly on $\theta$. Thus \eqref{eq:k} yields
\begin{equation*}
k_u\in C^{2,\delta}(\mathbb{S}^1),\quad k_u>0\quad\text{on }\mathbb{S}^1.
\end{equation*}
Using $\mathrm{d}\theta=\kappa\mathrm{d}s$ in \eqref{eq:flux}, we have
\begin{equation*}
\int_0^{2\pi}k_u(\theta)n(\theta)\mathrm{d}\theta=\frac{1}{2}\int_{\partial\Omega}q_u^2\nu\mathrm{d}s=0.
\end{equation*}

By Proposition \ref{prop:real}, there exists $h\in\mathcal{A}$ such that $\mathcal{K}(h)=k_u$. Set $\Omega_2=\Omega_h$ and let $v=v_h$ be the positive ground state of \eqref{eq:le} in $\Omega_2$. The domain $\Omega_2$ has $C^{3,\delta}$ boundary of positive curvature, and \eqref{eq:map} gives $k_v=k_u$ on $\mathbb{S}^1$. Lemma \ref{lem:conv} implies $D^2(v^{-\frac{p-1}{2}})>0$ in $\Omega_2$. By Proposition \ref{prop:comp}, there exists $c\in\mathbb{R}^2$ such that
\begin{equation*}
\Omega=c+\Omega_2,\quad u(x)=v(x-c)\quad\text{for all }x\in\Omega.
\end{equation*}
Since $\Omega=c+\Omega_2$, we have $S_p(\Omega)=S_p(\Omega_2)$ by \eqref{eq:sob}. Changing variables $y=x-c$ and using \eqref{eq:least} for $v$, we obtain
\begin{equation*}
\int_\Omega u^{p+1}\mathrm{d}x=\int_{\Omega_2}v^{p+1}\mathrm{d}y=S_p(\Omega_2)^{\frac{p+1}{p-1}}=S_p(\Omega)^{\frac{p+1}{p-1}}.
\end{equation*}
Therefore, $u$ is a ground state of problem \eqref{eq:le} in $\Omega$.
\end{proof}

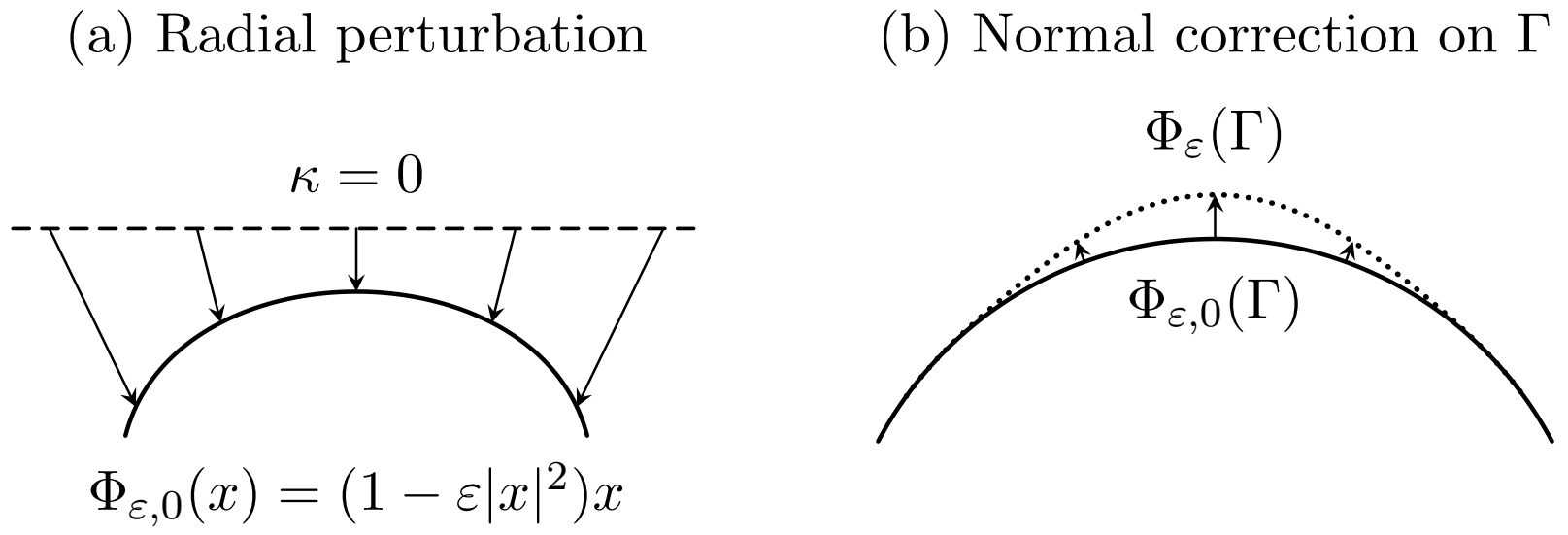
\begin{figure}[htbp]
\centering
\begin{tikzpicture}[
font=\small, line cap=round, line join=round, >=stealth,
original/.style={dash pattern=on 3.2pt off 2.4pt,line width=0.7pt},
radial/.style={solid,line width=0.85pt},
correction/.style={dash pattern=on 0pt off 2.25pt,line width=1.05pt},
motion/.style={->,line width=0.48pt}
]
\def\radialfactor{(1-0.10*(1+(\t)*(\t)))}
\def\correctedradius{(2.6+0.30*(cos(90*(\t)/54))^4)}
\begin{scope}
\node at (0,1.85) {(a) Radial perturbation};
\draw[original] (-2.34,0.55) -- (2.34,0.55);
\node[above=3pt] at (0,0.55) {$\kappa=0$};
\foreach \t in {-1.35,-0.7,0,0.7,1.35}{
\draw[motion] ({1.55*\t},0.55)
-- ({1.55*\t*\radialfactor},{0.55+4.3*(\radialfactor-1)});
}
\draw[radial] plot[domain=-1.51:1.51,samples=91,variable=\t]
({1.55*\t*\radialfactor},{0.55+4.3*(\radialfactor-1)});
\node at (0,-1.26) {$\Phi_{\varepsilon,0}(x)=(1-\varepsilon|x|^2)x$};
\end{scope}

\begin{scope}[xshift=5.85cm]
\node at (0,1.85) {(b) Normal correction on $\Gamma$};
\draw[radial] plot[domain=-62:62,samples=125,variable=\t]
({2.6*sin(\t)},{2.6*cos(\t)-2.6+0.48});
\foreach \t in {-20,0,20}{
\draw[motion] ({2.6*sin(\t)},{2.6*cos(\t)-2.6+0.48})
-- ({\correctedradius*sin(\t)},{\correctedradius*cos(\t)-2.6+0.48});
}
\draw[correction] plot[domain=-54:54,samples=109,variable=\t]
({\correctedradius*sin(\t)},{\correctedradius*cos(\t)-2.6+0.48});
\node[below=3pt] at (0,0.48) {$\Phi_{\varepsilon,0}(\Gamma)$};
\node[above=2pt] at (0,0.78) {$\Phi_\varepsilon(\Gamma)$};
\end{scope}
\end{tikzpicture}
\caption{The two boundary perturbations.}
\label{fig:domain-per}
\end{figure}

When the boundary curvature vanishes somewhere, we first perturb $\Omega$ to make the curvature positive. If $\ker L_u\neq\{0\}$, we also vary the boundary on an arc where $\kappa>0$ to approximate the given solution $u$ of problem \eqref{eq:le}. Figure \ref{fig:domain-per} illustrates the two perturbations.

\begin{lemma}\label{lem:approx}
Let $u$ be a positive solution of problem \eqref{eq:le} in a bounded smooth convex domain $\Omega\subset\mathbb{R}^2$, and fix $0<\delta<\min\{1,p-1\}$. For every sufficiently small $\varepsilon>0$, there exists a smooth diffeomorphism $\Phi_\varepsilon$ of $\mathbb{R}^2$ such that $\Omega_\varepsilon=\Phi_\varepsilon(\Omega)$ is convex and has positive boundary curvature. Moreover, problem \eqref{eq:le} has a positive solution $u_\varepsilon$ in $\Omega_\varepsilon$ satisfying
\begin{equation}\label{eq:approx}
\|\Phi_\varepsilon-\mathrm{Id}\|_{C^2(\mathbb{R}^2)}+\|u_\varepsilon\circ\Phi_\varepsilon-u\|_{C^{2,\delta}(\overline{\Omega})}\to0\quad\text{as }\varepsilon\to0^+.
\end{equation}
\end{lemma}

\begin{proof}
After a translation, we assume that $0\in\Omega$. Let $x(s)$ parametrize $\partial\Omega$ by arc length in the direction of $\tau$, with period $L=|\partial\Omega|$. Since $x'(s)=\tau(x(s))$, \eqref{eq:Fre} gives
\begin{equation*}
x''(s)=(\partial_\tau\tau)(x(s))=-\kappa(x(s))\nu(x(s)).
\end{equation*}
For each $s_0$, using the convexity of $\Omega$, the function $s\mapsto x(s)\cdot\nu(x(s_0))$ attains its maximum at $s_0$. Hence
\begin{equation*}
\kappa(x(s))=-x''(s)\cdot\nu(x(s))\geq0\quad\text{for all }s.
\end{equation*}
If $\kappa\equiv0$, then $x'(s)=x'(0)$ for all $s$. By periodicity,
\begin{equation*}
0=x(L)-x(0)=Lx'(0),
\end{equation*}
contrary to $|x'(0)|=1$. Thus $\kappa$ is positive at some point of $\partial\Omega$. By continuity of $\kappa$, there exists an open arc $\Gamma\subset\partial\Omega$ and a constant $\kappa_0>0$ such that
\begin{equation*}
\kappa(x)\geq\kappa_0\quad\text{for all }x\in\overline{\Gamma}.
\end{equation*}

Choose $V\in C_c^\infty(\mathbb{R}^2;\mathbb{R}^2)$ with $V(x)=-|x|^2x$ near $\overline{\Omega}$. Set
\begin{equation}\label{eq:var}
\Phi_{\varepsilon,0}=\mathrm{Id}+\varepsilon V,\quad x_\varepsilon(s)=\Phi_{\varepsilon,0}(x(s)).
\end{equation}
If $|\varepsilon|\|DV\|_{L^\infty}<1$, the contraction mapping theorem gives a unique solution of $x=y-\varepsilon V(x)$ for each $y\in\mathbb{R}^2$. Thus $\Phi_{\varepsilon,0}$ is bijective. Since $I+\varepsilon DV$ is invertible at every point, $\Phi_{\varepsilon,0}$ is a smooth diffeomorphism of $\mathbb{R}^2$.

\vskip 0.1cm
Differentiating \eqref{eq:var} with respect to $s$ and using $x'=\tau$ and $x''=-\kappa\nu$, we obtain
\begin{equation*}
x_\varepsilon'=(I+\varepsilon DV)\tau,\quad
x_\varepsilon''=-\kappa(I+\varepsilon DV)\nu+\varepsilon D^2V[\tau,\tau].
\end{equation*}
Since $V(x)=-|x|^2x$ near $\partial\Omega$ and $\det(\tau,x)=-x\cdot\nu$, we have
\begin{equation*}
\begin{aligned}
D^2V[\tau,\tau]&=-2x-4(x\cdot\tau)\tau,~~~~~~~
\det(\tau,D^2V[\tau,\tau]) =-2\det(\tau,x)=2x\cdot\nu.
\end{aligned}
\end{equation*}
Using $\det(\tau,\nu)=-1$, we find
\begin{equation*}
\begin{aligned}
\det(x_\varepsilon',x_\varepsilon'')&=\kappa\det(I+\varepsilon DV)+\varepsilon\det((I+\varepsilon DV)\tau,D^2V[\tau,\tau])\\
&=\kappa\det(I+\varepsilon DV)+2\varepsilon x\cdot\nu+O(\varepsilon^2).
\end{aligned}
\end{equation*}
Set $d=\operatorname{dist}(0,\partial\Omega)$. As in the proof of Lemma \ref{lem:tensor}, convexity gives $x\cdot\nu(x)\geq d>0$ on $\partial\Omega$. Since $I+\varepsilon DV\to I$ uniformly, for all sufficiently small $\varepsilon>0$, we have
\begin{equation*}
\det(x_\varepsilon',x_\varepsilon'')\geq\frac{\kappa(x(s))}{2}+\varepsilon d,\quad \frac{1}{2}\leq|x_\varepsilon'|\leq2.
\end{equation*}
By the curvature formula \cite[Section 1-5]{docarmo}, we obtain
\begin{equation}\label{eq:curv}
\kappa_{\varepsilon,0}(x_\varepsilon(s))=\frac{\det(x_\varepsilon',x_\varepsilon'')}{|x_\varepsilon'|^3}\geq c\left(\kappa(x(s))+\varepsilon\right), \quad c=\frac{1}{8}\min\left\{\frac{1}{2},d\right\}.
\end{equation}

Set $E=C^{2,\delta}_0(\overline{\Omega})$ and $\ell=\dim\ker L_u$, with $L_u$ acting on $E$. Since $pu^{p-1}\in C^{0,\delta}(\overline{\Omega})$, the argument in the proof of Lemma \ref{lem:shape} shows that $L_u(E)$ is closed and has codimension $\ell$ in $C^{0,\delta}(\overline{\Omega})$.
Choose an $L^2(\Omega)$-orthonormal basis $e_1,\cdots,e_\ell$ of $\ker L_u$. For $z\in E$, Green's identity yields
\begin{equation*}
\int_\Omega (L_uz)e_i\mathrm{d}x=\int_\Omega zL_ue_i\mathrm{d}x=0,\quad 1\leq i\leq\ell.
\end{equation*}
As $L_u(E)$ has codimension $\ell$, for every $f\in C^{0,\delta}(\overline{\Omega})$, we have
\begin{equation*}
f\in L_u(E)\quad\Longleftrightarrow\quad \int_\Omega fe_i\mathrm{d}x=0\quad\text{for }1\leq i\leq\ell.
\end{equation*}
Let $E_0$ be the subspace of $E$ orthogonal to $\ker L_u$ in $L^2(\Omega)$. Since $E=E_0\oplus\ker L_u$, the restriction of $L_u$ to $E_0$ is an isomorphism onto $L_u(E)$.

\vskip 0.1cm
Suppose $\ell>0$. We seek $g_1,\cdots,g_\ell\in C_c^\infty(\Gamma)$ satisfying
\begin{equation}\label{eq:orth}
\int_\Gamma g_j\partial_\nu u\partial_\nu e_i\mathrm{d}s=\delta_{ij},\quad 1\leq i,j\leq\ell.
\end{equation}
For $g\in C_c^\infty(\Gamma)$, set
\begin{equation*}
\mathcal{B}(g)=\left(\int_\Gamma g\partial_\nu u\partial_\nu e_i\mathrm{d}s\right)_{i=1}^\ell.
\end{equation*}
Let $a=(a_1,\cdots,a_\ell)\in\mathbb{R}^\ell$ be orthogonal to $\operatorname{Range}\mathcal{B}$. Setting $e=\displaystyle\sum_{i=1}^\ell a_ie_i$, we obtain
\begin{equation*}
0=a\cdot\mathcal{B}(g)=\int_\Gamma g\partial_\nu u\partial_\nu e\mathrm{d}s\quad\text{for every }g\in C_c^\infty(\Gamma).
\end{equation*}
Since $\partial_\nu u<0$ on $\Gamma$ by \eqref{eq:reg}, we have $\partial_\nu e=0$ on $\Gamma$. As $e\in\ker L_u$,
\begin{equation*}
\begin{cases}
-\Delta e=pu^{p-1}e &\text{in }\Omega,\\
e=\partial_\nu e=0 &\text{on }\Gamma.
\end{cases}
\end{equation*}
The unique continuation principle implies $e=0$ in $\Omega$. Hence $\operatorname{Range}\mathcal{B}=\mathbb{R}^\ell$ and \eqref{eq:orth} holds.

\vskip 0.1cm

Since $g_j \in C_c^\infty(\Gamma)$, we extend each $g_j$ by zero to $\partial\Omega$. As in the proof of Lemma \ref{lem:shape},
\begin{equation*}
x=\xi+t\nu(\xi),\quad \xi\in\partial\Omega,\quad |t|<t_0,
\end{equation*}
where $t_0>0$ is sufficiently small. Choose $\chi\in C_c^\infty((-t_0,t_0))$ with $\chi(0)=1$ and set
\begin{equation*}
V_j(\xi+t\nu(\xi))=\chi(t)g_j(\xi)\nu(\xi).
\end{equation*}
Extend $V_j$ by zero outside this neighborhood. Then $V_j\in C_c^\infty(\mathbb{R}^2;\mathbb{R}^2)$ and $V_j=g_j\nu$ on $\partial\Omega$. Since $g_j$ has compact support in $\Gamma$, $V_j$ vanishes near $\partial\Omega\setminus\Gamma$.

\vskip 0.1cm
For $b=(b_1,\cdots,b_\ell)\in\mathbb{R}^\ell$, set
\begin{equation*}
\Phi_{\varepsilon,b}=\Phi_{\varepsilon,0}+\sum_{j=1}^\ell b_jV_j.
\end{equation*}
For small $|\varepsilon|+|b|$, the argument for $\Phi_{\varepsilon,0}$ in \eqref{eq:var} then shows that $\Phi_{\varepsilon,b}$ is a smooth diffeomorphism of $\mathbb{R}^2$. When $\ell=0$, the sum is empty and $E_0=E$.
For $(\varepsilon,b,z)\in\mathbb{R}\times\mathbb{R}^\ell\times E_0$ near $(0,0,0)$, set
\begin{equation} \label{eq:F}
\mathcal{F}(\varepsilon,b,z)=\left[-\Delta\left((u+z)\circ\Phi_{\varepsilon,b}^{-1}\right)\right]\circ\Phi_{\varepsilon,b}-F(u+z).
\end{equation}
By \eqref{eq:reg}, Lemma \ref{lem:quot} applies with $v=u$. Hence
\begin{equation}\label{eq:quotu}
\left\|\frac{z}{u}\right\|_{C^{1,\delta}(\overline{\Omega})}\leq C\|z\|_{C^{2,\delta}(\overline{\Omega})}\quad\text{for every }z\in E.
\end{equation}
By Lemma \ref{lem:quot}, $z\mapsto F(u+z)$ is smooth near zero in $E$ with values in $C^{0,\delta}(\overline{\Omega})$. Together with \eqref{eq:pull}, this yields the smoothness of $\mathcal{F}$ near $(0,0,0)$. At $(0,0,0)$, we have
\begin{equation*}
\mathcal{F}(0,0,0)=0,\quad D_z\mathcal{F}(0,0,0)[\zeta]=L_u\zeta\quad\text{for }\zeta\in E_0.
\end{equation*}

For $1\leq j\leq\ell$, set $W_j=\partial_{b_j}\mathcal{F}(0,0,0)$. By \eqref{eq:opder}, we have
\begin{equation*}
W_j=-L_u(V_j\cdot\nabla u).
\end{equation*}
Set $\varphi_j=V_j\cdot\nabla u$. Since $W_j=-L_u\varphi_j$ and $\Delta e_i=-pu^{p-1}e_i$, we have
\begin{equation*}
e_iW_j=e_i\Delta\varphi_j+pu^{p-1}e_i\varphi_j=e_i\Delta\varphi_j-\varphi_j\Delta e_i.
\end{equation*}
Applying Green's identity with $e_i=0$ on $\partial\Omega$, we obtain
\begin{equation*}
\int_\Omega e_iW_j\mathrm{d}x=\int_{\partial\Omega}\left(e_i\partial_\nu\varphi_j-\varphi_j\partial_\nu e_i\right)\mathrm{d}s=-\int_{\partial\Omega}(V_j\cdot\nabla u)\partial_\nu e_i\mathrm{d}s.
\end{equation*}
Since $V_j=g_j\nu$ on $\partial\Omega$ and $g_j\in C_c^\infty(\Gamma)$ satisfying \eqref{eq:orth},
\begin{equation*}
\int_\Omega e_iW_j\mathrm{d}x=-\int_\Gamma g_j\partial_\nu u\partial_\nu e_i\mathrm{d}s=-\delta_{ij}.
\end{equation*}
Therefore, for $f\in C^{0,\delta}(\overline{\Omega})$, the equation $D_{(b,z)}\mathcal{F}(0,0,0)[b,z]=f$ is equivalent to
\begin{equation*}
L_uz=f-\sum_{j=1}^\ell b_jW_j,\quad z\in E_0.
\end{equation*}
The right-hand side belongs to $L_u(E)$ if and only if
\begin{equation*}
0=\int_\Omega\left(f-\sum_{j=1}^\ell b_jW_j\right)e_i\mathrm{d}x=\int_\Omega fe_i\mathrm{d}x+b_i,\quad 1\leq i\leq\ell.
\end{equation*}
Thus $b$ is uniquely determined by
\begin{equation*}
b_i=-\int_\Omega fe_i\mathrm{d}x,\quad 1\leq i\leq\ell.
\end{equation*}
For this $b$, the isomorphism $L_u:E_0\to L_u(E)$ yields a unique solution $z\in E_0$. Hence $D_{(b,z)}\mathcal{F}(0,0,0)$ is an isomorphism.
For $|\varepsilon|$ sufficiently small, the implicit function theorem yields a unique pair $(b_\varepsilon,z_\varepsilon)$ near $(0,0)$ satisfying
$
\mathcal{F}(\varepsilon,b_\varepsilon,z_\varepsilon)=0.
$
The map $\varepsilon\mapsto(b_\varepsilon,z_\varepsilon)$ is $C^1$ and $(b_0,z_0)=(0,0)$. Hence
\begin{equation}\label{eq:ift}
|b_\varepsilon|+\|z_\varepsilon\|_{C^{2,\delta}(\overline{\Omega})}\leq C|\varepsilon|.
\end{equation}

Set $\Phi_\varepsilon=\Phi_{\varepsilon,b_\varepsilon}$ and $\Omega_\varepsilon=\Phi_\varepsilon(\Omega)$. By \eqref{eq:var}, the definition of $\Phi_{\varepsilon,b}$ and \eqref{eq:ift},
\begin{equation} \label{eq:Phi-Id}
\|\Phi_\varepsilon-\mathrm{Id}\|_{C^2(\mathbb{R}^2)}\leq C\left(|\varepsilon|+|b_\varepsilon|\right)\leq C|\varepsilon|.
\end{equation}
Let $y_\varepsilon(s)=\Phi_\varepsilon(x(s))$. Using $x'=\tau$ and $x''=-\kappa\nu$ from \eqref{eq:Fre}, we obtain
\begin{equation*}
\begin{aligned}
y_\varepsilon'-x'&=\left(D\Phi_\varepsilon(x(s))-I\right)\tau,~~~~~~
y_\varepsilon''-x'' =D^2\Phi_\varepsilon(x(s))[\tau,\tau]-\kappa\left(D\Phi_\varepsilon(x(s))-I\right)\nu.
\end{aligned}
\end{equation*}
Since $\kappa$ is bounded on $\partial\Omega$,
\begin{equation*}
\|y_\varepsilon-x\|_{C^2([0,L])}\leq C\|\Phi_\varepsilon-\mathrm{Id}\|_{C^2(\mathbb{R}^2)}\leq C|\varepsilon|.
\end{equation*}
For $x(s)\in\Gamma$, we have $\det(x',x'')=\kappa(x(s))\geq\kappa_0$. Hence, for all sufficiently small $\varepsilon>0$,
\begin{equation*}
\begin{aligned}
\det(y_\varepsilon',y_\varepsilon'')
&=\det(x',x'')+\det(y_\varepsilon'-x',y_\varepsilon'')+\det(x',y_\varepsilon''-x'')\\
&\geq\det(x',x'')-|y_\varepsilon'-x'||y_\varepsilon''|-|x'||y_\varepsilon''-x''|\\
&\geq\kappa_0-C\varepsilon\geq\frac{\kappa_0}{2}.
\end{aligned}
\end{equation*}
Thus the curvature is positive on $\Phi_\varepsilon(\Gamma)$. If $x(s)\in\partial\Omega\setminus\Gamma$, all $V_j$ vanish near $x(s)$, so $y_\varepsilon=x_\varepsilon$ near $s$. By \eqref{eq:curv},
\begin{equation*}
\frac{\det(y_\varepsilon',y_\varepsilon'')}{|y_\varepsilon'|^3}=\kappa_{\varepsilon,0}(x_\varepsilon(s))\geq c\left(\kappa(x(s))+\varepsilon\right)>0.
\end{equation*}
Since $\Phi_\varepsilon$ is a diffeomorphism, $y_\varepsilon=\Phi_\varepsilon\circ x$ parametrizes $\partial\Omega_\varepsilon=\Phi_\varepsilon(\partial\Omega)$ as a simple regular closed curve. Its curvature is positive everywhere, so $\Omega_\varepsilon$ is convex by \cite[Section 5-7]{docarmo}.

\vskip 0.1cm
Set $u_\varepsilon=(u+z_\varepsilon)\circ\Phi_\varepsilon^{-1}$ in $\Omega_\varepsilon$. By \eqref{eq:quotu} and \eqref{eq:ift}, for all sufficiently small $\varepsilon>0$,
\begin{equation*}
\left\|\frac{z_\varepsilon}{u}\right\|_{C(\overline{\Omega})}
\leq C\|z_\varepsilon\|_{C^{2,\delta}(\overline{\Omega})}
\leq C\varepsilon\leq\frac{1}{2}.
\end{equation*}
For each $y\in\Omega_\varepsilon$, denote $x=\Phi_\varepsilon^{-1}(y)\in\Omega$. Then
\begin{equation*}
u_\varepsilon(y)=u(x)\left(1+\frac{z_\varepsilon(x)}{u(x)}\right)
\geq\frac{u(x)}{2}>0.
\end{equation*}
Since $\mathcal{F}(\varepsilon,b_\varepsilon,z_\varepsilon)=0$, \eqref{eq:F} gives
\begin{equation*}
(-\Delta u_\varepsilon)(\Phi_\varepsilon(x))
=F(u(x)+z_\varepsilon(x))=(u(x)+z_\varepsilon(x))^p.
\end{equation*}
Using \eqref{eq:le}, $u=0$ on $\partial\Omega$. Since $z_\varepsilon\in E_0$, we also have $z_\varepsilon=0$ on $\partial\Omega$. Thus
\begin{equation*}
\begin{cases}
-\Delta u_\varepsilon=u_\varepsilon^p,\quad u_\varepsilon>0 &\text{in }\Omega_\varepsilon,\\
u_\varepsilon=0 &\text{on }\partial\Omega_\varepsilon.
\end{cases}
\end{equation*}
Finally, with $u_\varepsilon\circ\Phi_\varepsilon-u=z_\varepsilon$, \eqref{eq:ift} and \eqref{eq:Phi-Id}, we have
\begin{equation*}
\|\Phi_\varepsilon-\mathrm{Id}\|_{C^2(\mathbb{R}^2)}
+\|u_\varepsilon\circ\Phi_\varepsilon-u\|_{C^{2,\delta}(\overline{\Omega})}
\leq C\varepsilon\to0.
\end{equation*}
This proves \eqref{eq:approx}.
\end{proof}

We now prove Theorem \ref{thm:main} and Corollary \ref{cor:conv}.

\begin{proof}[{\bf{Proof of Theorem \ref{thm:main}}}]
Let $u$ be a positive solution of problem \eqref{eq:le} in $\Omega$, and let $\Omega_\varepsilon=\Phi_\varepsilon(\Omega)$ and $u_\varepsilon$ be as in Lemma \ref{lem:approx}. Since $\partial\Omega_\varepsilon$ has positive curvature, Proposition \ref{prop:pos} and \eqref{eq:least} imply
\begin{equation}\label{eq:limit}
\big(S_p(\Omega_\varepsilon)\big)^{\frac{p+1}{p-1}}
=\int_{\Omega_\varepsilon}u_\varepsilon^{p+1}\mathrm{d}y
=\int_\Omega (u_\varepsilon\circ\Phi_\varepsilon)^{p+1}|\det D\Phi_\varepsilon|\mathrm{d}x.
\end{equation}
By \eqref{eq:approx}, we get
\begin{equation*}
u_\varepsilon\circ\Phi_\varepsilon\to u,\quad \det D\Phi_\varepsilon\to1
\quad\text{uniformly on }\overline\Omega.
\end{equation*}
Since $S_p(\Omega_\varepsilon)\to S_p(\Omega)$ by \eqref{eq:scont}, passing to the limit in \eqref{eq:limit} gives
\begin{equation*}
\int_\Omega u^{p+1}\mathrm{d}x=\big(S_p(\Omega)\big)^{\frac{p+1}{p-1}}.
\end{equation*}
By \eqref{eq:least}, $u$ is a ground state of problem \eqref{eq:le} in $\Omega$. Since $u$ was an arbitrary positive solution, Proposition \ref{prop:bf} shows that problem \eqref{eq:le} has exactly one positive solution in $\Omega$. The inequalities $\lambda_1(L_u)<0<\lambda_2(L_u)$ follow from \eqref{eq:spec}.
\end{proof}

\begin{proof}[{\bf{Proof of Corollary \ref{cor:conv}}}]
By Theorem \ref{thm:main}, $u$ is the positive ground state of problem \eqref{eq:le} in $\Omega$. Applying Lemma \ref{lem:conv} with $\beta=\frac{p-1}{2}$, we obtain
\begin{equation*}
D^2\left(u^{-\frac{p-1}{2}}\right)>0\quad\text{in }\Omega.
\end{equation*}
\end{proof}

\medskip

\subsection*{Acknowledgements}\

\vskip 0.1cm

This work was supported by the National Key R\&D Program (No. 2023YFA1010002).
Hong-Ge Chen was supported by the National Natural Science Foundation of China (Grant Nos. 12671239 and 12571249).
Peng Luo and Shuangjie Peng were supported by the National Natural Science Foundation of China (Grant No. 12631005).

\vskip 0.2truein
\paragraph{\textbf{Declaration of AI-assisted tools.}}

\vskip 0.1cm
During the preparation of this manuscript, the authors used AI tools as auxiliary aids for exploratory discussions, preliminary proof drafting, language editing, and drawing a schematic figure. The authors take full responsibility for the mathematical validation, the final proof choices, and the final text.


\begin{thebibliography}{99}

\bibitem{ag}
Adimurthi and M. Grossi,
{\it Asymptotic estimates for a two-dimensional problem with polynomial nonlinearity},
Proc. Amer. Math. Soc. {\bf 132} (2004), no. 4, 1013--1019.

\bibitem{bf}
L. Brasco and G. Franzina,
{\it An overview on constrained critical points of Dirichlet integrals},
Rend. Semin. Mat. Univ. Politec. Torino {\bf 78} (2020), no. 2, 7--50.

\bibitem{brezis}
H. R. Brezis,
{\it Functional analysis, Sobolev spaces and partial differential equations},
Universitext, Springer, New York, 2011.

\bibitem{cfr}
T. F. Carroll, M. M. Fall and J. Ratzkin,
{\it On the rate of change of the sharp constant in the Sobolev-Poincar\'e{} inequality},
Math. Nachr. {\bf 290} (2017), no. 14-15, 2185--2197.

\bibitem{docarmo}
M. P. do Carmo,
{\it Differential geometry of curves \& surfaces},
revised \& updated second edition, Dover, Mineola, NY, 2016.

\bibitem{dgp}
L. Damascelli, M. Grossi and F. Pacella,
{\it Qualitative properties of positive solutions of semilinear elliptic equations in symmetric domains via the maximum principle},
Ann. Inst. H. Poincar\'e{} Anal. Non Lin\'eaire {\bf 16} (1999), no. 5, 631--652.

\bibitem{dancer}
E. N. Dancer,
{\it The effect of domain shape on the number of positive solutions of certain nonlinear equations},
J. Differential Equations {\bf 74} (1988), no. 1, 120--156.


\bibitem{dm}
F. De Marchis, M. Grossi, I. Ianni and F. Pacella,
{\it Morse index and uniqueness of positive solutions of the Lane-Emden problem in planar domains},
J. Math. Pures Appl. (9) {\bf 128} (2019), 339--378..

\bibitem{dip}
F. De Marchis, I. Ianni and F. Pacella,
{\it Asymptotic profile of positive solutions of Lane-Emden problems in dimension two},
J. Fixed Point Theory Appl. {\bf 19} (2017), no. 1, 889--916.

\bibitem{gms}
M. Gallo, S. J. N. Mosconi and M. Squassina,
{\it Power law convergence and concavity for the logarithmic Schr\"odinger equation},
Math. Ann. {\bf 395} (2026), no. 1, Paper No. 21, 60 pp.

\bibitem{gnn}
B. Gidas, W.-M. Ni and L. Nirenberg,
{\it Symmetry and related properties via the maximum principle},
Comm. Math. Phys. {\bf 68} (1979), no. 3, 209--243.

\bibitem{gt}
D. Gilbarg and N. S. Trudinger,
{\it Elliptic partial differential equations of second order}, reprint of the 1998 edition,
Classics in Mathematics, Springer, Berlin, 2001.

\bibitem{Gr}
P. Grisvard, \emph{Elliptic Problems in Nonsmooth Domains},
{\it Monographs and Studies in Mathematics, vol.~24},
Pitman (Advanced Publishing Program), Boston, MA, 1985.

\bibitem{gily}
M. Grossi, I. Ianni, P. Luo and S. Yan,
{\it Non-degeneracy and local uniqueness of positive solutions to the Lane-Emden problem in dimension two},
J. Math. Pures Appl. (9) {\bf 157} (2022), 145--210.

\bibitem{gtaka}
M. Grossi and F. Takahashi,
{\it Nonexistence of multi-bubble solutions to some elliptic equations on convex domains},
J. Funct. Anal. {\bf 259} (2010), no. 4, 904--917.


\bibitem{ks}
N. Kamburov and B. S. Sirakov,
{\it Uniform a priori estimates for positive solutions of the Lane-Emden equation in the plane},
Calc. Var. Partial Differential Equations {\bf 57} (2018), no. 6, Paper No. 164, 8 pp.

\bibitem{lin}
C.-S. Lin,
{\it Uniqueness of least energy solutions to a semilinear elliptic equation in ${\bf R}^2$},
Manuscripta Math. {\bf 84} (1994), no. 1, 13--19.

\bibitem{mardare}
S. Mardare,
{\it On Poincar\'e{} and de Rham's theorems},
Rev. Roumaine Math. Pures Appl. {\bf 53} (2008), no. 5-6, 523--541.

\bibitem{mira}
P. Mira,
{\it Overdetermined elliptic problems in topological disks},
J. Differential Equations {\bf 264} (2018), no. 11, 6994--7005.

\bibitem{pohoz}
S. I. Pohozaev,
{\it On the eigenfunctions of the equation $\Delta u+\lambda f(u)=0$},
Dokl. Akad. Nauk SSSR {\bf 165} (1965), 36--39.

\bibitem{pucci}
P. Pucci and J. B. Serrin Jr.,
{\it A general variational identity},
Indiana Univ. Math. J. {\bf 35} (1986), no. 3, 681--703.

\bibitem{ren94}
X. Ren and J. Wei,
{\it On a two-dimensional elliptic problem with large exponent in nonlinearity},
Trans. Amer. Math. Soc. {\bf 343} (1994), no. 2, 749--76.

\bibitem{ren96}
X. Ren and J. Wei,
{\it Single-point condensation and least-energy solutions},
Proc. Amer. Math. Soc. {\bf 124} (1996), no. 1, 111--120.

\bibitem{schneider}
R. Schneider,
{\it Convex bodies: the Brunn-Minkowski theory}, second expanded edition,
Encyclopedia of Mathematics and its Applications, 151, Cambridge Univ. Press, Cambridge, 2014.

\bibitem{simader}
C. G. Simader,
{\it On Dirichlet's boundary value problem},
Lecture Notes in Mathematics, Vol. 268, Springer, Berlin-New York, 1972.

\bibitem{dancer89}
E. N. Dancer,
{\it On the influence of domain shape on the existence of large solutions of some superlinear problems},
Math. Ann. {\bf 285} (1989), no. 4, 647--669.


\bibitem{lions}
P.-L. Lions,
{\it On the existence of positive solutions of semilinear elliptic equations},
SIAM Rev. {\bf 24} (1982), no. 4, 441--467.

\bibitem{lwz}
H. Li, J. Wei and W. Zou,
{\it Uniqueness, multiplicity and nondegeneracy of positive solutions to the Lane-Emden problem},
J. Math. Pures Appl. (9) {\bf 179} (2023), 1--67.

\end{thebibliography}
\end{document}